\documentclass{elsarticle}[preprint,10pt]

\journal{Journal of Computational Physics}

\usepackage{graphicx} 
\usepackage{algorithm}
\usepackage{algorithmic}

\usepackage[utf8]{inputenc}
\usepackage{subcaption}
\usepackage{graphicx}
\usepackage{bm}
\usepackage{comment}
\usepackage{psfrag}
\usepackage{latexsym,amsmath,amsfonts,amscd,amsthm}
\usepackage{changebar}
\usepackage{color}
\usepackage{bm}
\usepackage{tikz}
\usepackage{multirow}
\usepackage{xcolor}
\usepackage{hyperref}
\usepackage{setspace}
\usepackage{amssymb}
\usepackage{mathrsfs}
\usepackage{caption}
\usepackage{subcaption}

\usepackage[title]{appendix}
\usepackage{ulem}
\usepackage{tikz}
\usetikzlibrary{arrows,backgrounds}

\newlength{\spse}
\newtheorem{thm}{Theorem}[section]

\newtheorem{rem}[thm]{Remark}

\newcommand{\BU}{\mathbf{U}}
\newcommand{\BV}{\mathbf{V}}

\newcommand{\BF}{\mathbf{F}}

\newcommand{\BA}{\mathbf{A}}
\newcommand{\BLambda}{\mathbf{\Lambda}}

\newcommand{\bx}{\mathbf{x}}

\newcommand{\bw}{\mathbf{w}}
\newcommand{\bc}{\mathbf{c}}
\newcommand{\bb}{\mathbf{b}}
\newcommand{\bv}{\boldsymbol{v}}

\newcommand{\bn}{\boldsymbol{n}}

\newcommand{\half}{\frac{1}{2}}
\newcommand{\brho}{\boldsymbol{\rho}}
\newcommand{\bg}{\boldsymbol{g}}
\newcommand{\br}{\boldsymbol{r}}

\newcommand{\bG}{\boldsymbol{G}}
\newcommand{\wbG}{\widetilde{\boldsymbol{G}}}
\newcommand{\bff}{\boldsymbol{f}}

\newcommand{\mfw}{\mathfrak{w}}

\newcommand{\BJ}{{\bf{J}} }
\newcommand{\BH}{{\bf{H}} }
\newcommand{\BM}{{\bf{M}} }
\newcommand{\BD}{{\bf{D}}}
\newcommand{\BSigma}{{\mathbf{\Sigma}}}
\newcommand{\bmu}{{\boldsymbol{\mu}}}

\newcommand{\pzc}{\textcolor{black}}

\newcommand{\epssisa}{\epsilon_{\textrm{SISA}}}
\newcommand{\epspod}{\epsilon_{\textrm{POD}}}
\newcommand{\epsromsad}{\epsilon_{\textrm{ROMSAD}}}

\title{Reduced Order Model Enhanced  Source Iteration with Synthetic Acceleration for Parametric Radiative Transfer Equation}

\begin{document}

\author{
Zhichao Peng
}

\ead{pengzhic@ust.hk}

\affiliation{organization={Department of Mathematics, The Hong Kong University of Science and Technology},
            addressline={Clear Water
Bay, Kowloon},
city={Hong Kong},
country={China}
}

\begin{abstract}
Applications such as uncertainty quantification, shape optimization, and optical tomography, require solving the radiative transfer equation (RTE)  many times for various parameters. Efficient solvers for RTE are highly desired.

Source Iteration with Synthetic Acceleration (SISA) is \pzc{a}  popular and successful iterative solver for RTE. 
Synthetic Acceleration (SA) acts as a preconditioning step to accelerate the convergence of Source Iteration (SI). After each source iteration, classical SA strategies introduce a correction to the macroscopic particle density by solving a low order approximation to a kinetic correction equation. For example, Diffusion Synthetic Acceleration (DSA) uses the diffusion limit. However, these strategies may become less effective when the underlying low order approximations are not accurate enough. Furthermore, they do not exploit low rank structures concerning the parameters of parametric problems.

To address these issues, we propose enhancing SISA with data-driven ROMs for the parametric problem and the corresponding kinetic correction equation. 
First, the ROM for the parametric problem can be utilized to obtain an improved initial guess. 
Second, the ROM for the kinetic correction equation can be utilized to design a low rank approximation to it. Unlike the diffusion limit, this ROM-based approximation builds on the kinetic description of the correction equation and leverages low rank structures concerning the parameters. We further introduce
a novel SA strategy called ROMSAD. ROMSAD initially adopts our ROM-based approximation to exploit its greater efficiency in the early stage of SISA, and then automatically switches to DSA to leverage its robustness in the later stage.
Additionally, we propose an approach to construct the ROM for the kinetic correction equation without directly solving it. 

Through a series of numerical tests, we \pzc{compare} the proposed methods \pzc{with SI-DSA and DSA preconditioned Krylov solver}.  Particularly, for a multiscale parametric pin-cell problem, ROMSAD achieves approximately $10$ times the acceleration compared to \pzc{SI-DSA and $4$ times acceleration compared to DSA preconditioned GMRES}.
\end{abstract}


\maketitle

\section{Introduction}
Radiative transfer equation (RTE) is a kinetic equation describing the behavior of particles (such as photons and neutrons) propagating through a background medium. It plays an important role in medical imaging \cite{arridge2009optical}, nuclear engineering \cite{pomraning1973equations}, astrophysics \cite{janka2007theory}, and remote sensing \cite{spurr2001linearized}.
In applications such as shape optimization, uncertainty quantification and optical tomorgraphy, RTE needs to be solved many times for various parameters, such as boundary conditions, material properties or geometric configurations. As a result, efficient solvers for parametric RTE are highly desired.

Source Iteration with Synthetic Acceleration (SISA) is one of the most successful iterative solvers for RTE, developed over decades and widely applied in various applications. Instead of providing an exhaustive literature review for SISA, we refer readers to the review paper \cite{Adams2002FastIM}. It is well known that, without Synthetic Acceleration (SA), Source Iteration (SI)  may converge slowly for scattering dominant (optically thick) or  multiscale problems \cite{Adams2002FastIM}. SA can be seen as a preconditioning step to accelerate the convergence of SI by introducing a correction to the macroscopic particle density (also known as the scalar flux) after each source iteration. If the density correction is obtained by solving an ideal kinetic correction equation, SI will converge with at most two source iterations. However, solving this kinetic correction equation is as expensive as solving the original problem. In practice, a low order approximation to the kinetic correction equation is solved in the correction step of SA. For example, Diffusion Synthetic Acceleration (DSA) \cite{kopp1963synthetic,alcouffe1977dittusion,adams1992diffusion,wareing1993new} adopts its diffusion limit, Quasi-Diffusion method \cite{gol1964quasi,anistratov1993nonlinear,olivier2023family} uses the variable Eddington factor, and S$2$SA \cite{lorence1989s} employs a low order discrete ordinates ($S_N$) approximation. \pzc{To accelerate the convergence of the iterative solver, Krylov method have been applied in \cite{guthrie1999gmres} and further combined with preconditioners based on incomplete LU decomposition \cite{patton2002application}, Transport Syntheteic Acceleration \cite{zika2000transport}, and DSA \cite{warsa2004krylov}. We refer readers to the reivew paper \cite{azmy2010advances} for the development of Krylov method for RTE.}
Despite the success of \pzc{classical} SA strategies \pzc{as preconditioners for SI or Krylov method}, they still have some limitations. Their effectiveness relies on the accuracy of the underlying empirical low order approximation to the kinetic correction equation. Moreover, low rank structures with respect to the parameters of parametric problems are not exploited.

To go beyond these limitations, we propose to enhance SISA  by utilizing data-driven reduced order models (ROMs) for parametric RTE and the corresponding kinetic correction equation. Before presenting our methods, we briefly review the basic ideas of data-driven ROMs and their recent developments in the context of RTE. Data-driven ROMs typically follow an offline-online decomposition framework. In the offline stage, a low dimensional linear space is constructed by exploring low rank structures in the training data, i.e. solutions for parameters in a training set. In the online stage, reconstruction and prediction can be efficiently done through an interpolation or a projection based on the low dimensional space constructed offline. In recent years, ROMs for steady state and time dependent RTEs have been actively developed to utilize low rank structures in space-time domain \cite{choi2021space}, angular space \cite{buchan2015pod,hughes2020discontinuous,peng2022reduced},  angular-time domain \cite{peng2022micro}, parametric problems \cite{tano2021affine,behne2022minimally,behne2023parametric}  and eigenvalue problems \cite{mcclarren2019calculating}. ROMs leveraging the variable Eddington factor are proposed in \cite{coale2023reduced,coale2023reduced2,coale2023variable}. Besides data-driven ROMs, low rank tensor or matrix decompositions have been applied to design low rank solvers for RTE in \cite{peng2020low,peng2021high,kusch2021robust,einkemmer2021asymptotic,kusch2022low,dominesey2022reduced,yin2023semi,einkemmer2024asymptotic}.

Now, we outline our strategies to leverage ROMs to accelerate the convergence of SISA. 
\begin{enumerate}
    \item The ROM for the parametric problem can be viewed as a surrogate solver for the parametric problem. It can efficiently provide an improved initial guess for SISA.
    \item The ROM for the correction equation can be utilized to design a new SA strategy, called ROMSA. 
    Instead of using empirical low order approximations to the kinetic correction equation like its diffuion limit, ROMSA employs an approximation based on a data-driven ROM for the kinetic correction equation. This approximation directly builds on the kinetic description of the correction equation and exploits low rank structures with respect to the parameters of the underlying parametric problem. Additionally, we propose an approach to construct the ROM for the kinetic correction equation without directly solving it. 
\end{enumerate}
In our numerical tests, we observe that ROMSA achieves more significant acceleration than DSA in the early stage of SISA. However, it may suffer from an efficiency reduction and becomes slower than DSA as iterations continue. The cause of this efficiency reduction is as follows. As source iterations continue, source terms in the kinetic correction equations for the later iterations may exhibit significant shape variations compared to those in the early stage of SISA. However, including information for source terms corresponding to the later iterations in our training data may lead to prohibitive memory costs. To improve the robustness of ROMSA without sacrificing its high efficiency in the early stage of SISA, we design a SA strategy called ROMSAD, which combines ROMSA and DSA. ROMSAD adopts ROMSA in the first few iterations of SISA to leverage its high efficiency in the early stage, and then automatically switches to DSA. In our numerical tests, we observe that, overall, ROMSAD is more efficient than DSA and more robust than ROMSA. Specifically, ROMSAD achieves approximately $10$ times the acceleration compared to DSA for a parametric multiscale pin-cell problem. 

To contextualize our methods, we briefly review other methods utilizing data-driven ROMs or similar ideas to accelerate iterative solvers for RTE. For nonlinear RTEs, Dynamic Mode Decomposition (DMD) is exploited as a low rank update strategy for the SI in \cite{mcclarren2022data}. A neural network surrogate for the transport sweep in SI is developed in \cite{tano2021sweep}. These  methods focus on SI, while we concentrate on the correction step of SA.  Random Singular Value Decomposition (RSVD) has been applied to build a low-rank boundary-to-boundary map of a Schwartz solver for RTE \cite{chen2021low}. A fast solver, applying offline-online decomposition but not built on data-driven ROMs, is proposed in   \cite{fu2024fast}. This method is under the framework of the Tailored Finite Point Method (TFPM). The offline stage of this method can be seen as building an efficient preconditioner for RTE based on a matrix factorization exploiting local structures given by TFPM. 

The rest of this paper is organized as follows. In Sec. \ref{sec:background}, we briefly review discrete ordinates angular discretization,  upwind discontinuous Galerkin spatial  discretization, and the SISA iterative solver for the steady state RTE. In Sec. \ref{sec:romig-romsa}, we build ROMs for parametric RTE and the corresponding correction equation, and introduce our ROM-based enhancement for SISA. In Sec. \ref{sec:numerical}, the performance of the proposed methods is demonstrated through a series of numerical tests. In  Sec. \ref{sec:conclusion}, we draw our conclusions and outline potential future directions.

\section{Background\label{sec:background}}
The steady state linear RTE with one energy group, isotropic scattering and isotropic inflow boundary conditions is:
\begin{subequations}
\label{eq:rte}
    \begin{align}
    &\bv \cdot \nabla_{\bx} f(\bx,\bv) + \sigma_t(\bx) f(\bx,\bv) = \sigma_s(\bx) \rho(\bx) + G(\bx), \quad \rho(\bx) = \frac{1}{4\pi}\int_{\bv\in\mathbb{S}^2} f d\bv, \quad \bx\in\Omega_{\bx},\\
    &f(\bx,\bv) = g(\bx), \quad \bx\in \partial\Omega_{\bx},\;\bv\cdot \bn(\bx)<0.
    \end{align}
\end{subequations}
Here, $f(\bx,\bv)$ is the particle distribution (also known as the angular flux) with angular direction $\bv\in\mathbb{S}^2$ at spatial location $\bx\in\Omega_{\bx}$, $\rho(\bx)=\frac{1}{4\pi}\int_{\bv\in\mathbb{S}^2} f(\bx,\bv) d\bv$ is the macroscopic density (also called scalar flux), $\sigma_s(\bx)\geq 0$ is the scattering cross section, $\sigma_t(\bx)\geq \sigma_s(\bx)\geq 0$ is the total cross section, $\sigma_a(\bx) = \sigma_t(\bx)-\sigma_s(\bx)$ is the absorption cross section, $G(\bx)$ is an isotropic source, and $\bn(\bx)$ is the outward normal direction of $\Omega_{\bx}$ at $\bx\in\partial\Omega_{\bx}$.

\textbf{1D slab geometry:} Under symmetry assumptions, RTE \eqref{eq:rte} can be further simplified in 1D slab geometry: 
\begin{equation}
\xi \partial_x f(x,\xi)+\sigma_t(x)f(x,\xi)=\sigma_s(x)\rho(x)+G(x),\quad
\rho(x)=\frac{1}{2}\int_{-1}^1 f(x,\xi)d\xi,
\label{eq:slab_geometry}
\end{equation}
where the particle distribution $f(x,\xi)$ depends on location $x\in\mathbb{R}$ and cosine of the angle between angular direction $\bv$ and the $x$-axis, namely $\xi\in[-1,1]$.

In this section, we briefly review discrete ordinates angular discretization, upwind discontinuous Galerkin (DG) spatial discretization, and  Source Iteration with Synthetic Acceleration (SISA).
\subsection{Discrete ordinates ($S_N$) angular discretization}
We apply discrete ordinates ($S_N$) method \cite{pomraning1973equations} in angular space. Let  $\{(\bv_j,\omega_j)\}_{j=1}^{N_{\bv}}$ be a set of quadrature points in angular space and the corresponding quadrature weights satisfying $\sum_{j=1}^{N_{\bv}}\omega_j=1$. 
RTE is solved at these quadrature points by approximating the normalized integral term, $\rho(\bx)=\frac{1}{4\pi}\int_{\bv\in\mathbb{S}^2}f(\bx,\bv)d\bv$, with the associated quadrature rule:
\begin{subequations}
\label{eq:sn_system}
\begin{align}
&(\bv_j\cdot\nabla_x+\sigma_t(\bx) )f(\bx,\bv_j) = \sigma_s(\bx) \rho(\bx)+G(\bx),\quad \rho(\bx)= \sum_{j=1}^{N_{\bv}} \omega_j f(\bx,\bv_j)\\
&f(\bx,\bv_j) = g(\bx), \quad \bx\in\partial\Omega_{\bx},\;\bv_j\cdot\bn(\bx)<0.
\end{align}
\end{subequations}

For 1D slab geometry, we use Gauss-Legendre quadrature points. For 2D X-Y geometry with $\mathbb{S}^2$ as the angular space, \pzc{we use Chebyshev-Legendre (CL) points}. CL quadrature is the tensor product of Chebyshev rule for the unit circle and Gauss-Legendre rule for $\bv_z\in[-1,1]$. The quadrature points and weights of the $N_{\phi}$-points normalized Chebyshev quadrature rule for the unit circle is 
\begin{equation*}
\left\{(\phi_j,\omega_{\phi,j}): \;\phi_j = \frac{2j\pi}{N_\phi}-\frac{\pi}{N_\phi}\;\text{and}\;\omega^\phi_j=\frac{1}{N_\phi},\; j=1,\dots,N_\phi\right\}.
\end{equation*}
Let  $\{(v_{z,j},\omega_{z,j})\}_{j=1}^{N_{v_z}}$ denote quadrature points and weights of the $N_{v_z}$-points normalized Gauss-Legendre rule for $[-1,1]$. The quadrature points and weights  of the CL$(N_\phi,N_{v_z})$ quadrature rule, $(\bv_j,\omega_j)$,   are defined as  
\begin{equation*}
\bv_j=\left(\cos(\phi_{j_1})\sqrt{1-v_{z,j_2}^2},\sin(\phi_{j_1})\sqrt{1-v_{z,j_2}^2},v_{z,j_2}\right),\;\omega_j=\omega_{\phi,j_1}\omega_{v_z,j_2}, 
\end{equation*}
where $1\leq j_1\leq N_\phi$, $1\leq j_2\leq N_{v_z}$ and $j=(j_2-1) N_{\phi}+j_1$.
Normalized integral on $\mathbb{S}^2$ can be approximated by this quadrature rule as 
$\frac{1}{4\pi}\int_{\bv\in\mathbb{S}^2}f(\bv)d\bv\approx \sum_{j=1}^{N_\phi N_{v_z}}\omega_j f(\bv_j).$

\subsection{Upwind discontinuous Galerkin spatial discretization}
We apply upwind discontinuous Galerkin (DG) spatial discretization, because it is  an asymptotic preserving (locking free in optically thick regions)  scheme \cite{adams2001discontinuous,guermond2010asymptotic}, which is able to capture the correct diffusion limit without resolving small mean free path of particles.

Consider 2D $X$-$Y$ geometry with a rectangular computational domain $\Omega_{\bx}=[x_l,x_r]\times[y_l,y_r]$. Let $\mathcal{T}_h=\{\mathcal{T}_i\}_{i=1}^{N_{\bx}}$ be a partition of $\Omega_\bx$ with $\mathcal{T}_i$'s being rectangles. We seek the solution in the discrete space
\begin{equation}
       U_h^K(\mathcal{T}_h):=\{u(\bx): u(\bx)|_{\mathcal{T}_{i}}\in Q^K(\mathcal{T}_{i}),1\leq i \leq N_{\bx}\},\label{eq:discrete_space}
\end{equation}
where  $Q^K(\mathcal{T}_{i})$ is the bi-variate polynomial space whose degree in each direction is
at most $K$ on the element $\mathcal{T}_{i}$. Denote the set of cell edges as $\partial \mathcal{T}_h$ and the set of edges on the inflow boundary for $\bv_j$ as 
\begin{equation*}
\partial \mathcal{T}_{h,j}^{(\textrm{ibc})}=\{\mathcal{E}: \mathcal{E}\in\partial\mathcal{T}_h,\;\mathcal{E}\subset\partial\Omega_{\bx},\;\bv_j\cdot\bn(\bx)<0,\forall \bx\in\mathcal{E}\},
\end{equation*}
where $\bn(\bx)$ is the outward normal direction of $\Omega_{\bx}$ at $\bx$.

Applying upwind DG spatial discretization to \eqref{eq:sn_system}, we seek $f_h(\bx,\bv_j)\in U_h^K(\mathcal{T}_h)$, $j=1,\dots,N_{\bv}$ satisfying $\forall\phi_h(\bx)\in U_h^K(\mathcal{T}_h),$
\begin{align}
-\sum_{i=1}^{N_{\bx}}&\int_{\mathcal{T}_i} \Big(\bv_j\cdot\nabla\phi_h(\bx)\Big) f_h(\bx,\bv_j) d\bx+\sum_{\mathcal{E}\in\partial\mathcal{T}_h\setminus\partial\mathcal{T}_{h,j}^{(\textrm{ibc})}} \int_{\mathcal{E}} \widehat{\BH}(\bv_j,f_h, \bn(\bx))\phi_h(\bx) d\bx+\sum_{i=1}^{N_{\bx}}\int_{\mathcal{T}_i}\sigma_t(\bx) f_h(\bx,\bv_j)\phi_h(\bx) d\bx
\notag
\\
= &\sum_{i=1}^{N_{\bx}}\int_{\mathcal{T}_i}\sigma_s(\bx) \rho_h(\bx)\phi_h(\bx) d\bx
 + \sum_{i=1}^{N_{\bx}}\int_{\mathcal{T}_i}G(\bx)\phi_h(\bx) d\bx
 -\sum_{\mathcal{E}\in\partial\mathcal{T}_{h,j}^{(\textrm{ibc})}}\int_{\mathcal{E}} g(\bx)\phi_h(\bx) \bv_j\cdot \bn(\bx) d\bx.
 \label{eq:DG}
\end{align}
Here, $ \rho_h(\bx) = \sum_{j=1}^{N_{\bv}}\omega_j f_h(\bx,\bv_j)$. The upwind numerical flux $\widehat{\BH}(\bv_j, f_h,\bn(\bx))$ along the edge $\mathcal{E}$ for an element $\mathcal{T}_i=\mathcal{T}^-$ with the neighboring element $\mathcal{T}^+$, is defined as
\begin{align}
\widehat{\BH}(\bv_j, f_h,\bn(\bx)) = \frac{\bv_j\cdot\bn(\bx)}{2}\Big(f_h^+(\bx,\bv_j)+f_h^-(\bx,\bv_j)\Big)+\frac{|\bv_j\cdot \bn(\bx)|}{2}\Big(f_h^-(\bx,\bv_j)-f_h^+(\bx,\bv_j)\Big),
\label{eq:upwind}
\end{align}
where $f_h^{\pm}(\bx,\bv_j)$ is the restriction of $f_h(\bx,\bv_j)$ to $\mathcal{T}^{\pm}$, and $\bn(\bx)$ is  the unit outward normal direction  at $\bx\in\mathcal{E}$ with respect to the element $\mathcal{T}_i=\mathcal{T}^-$. 

We further rewrite the DG scheme to its matrix-vector form. Let $\{\phi_i(\bx)\}_{i=1}^{N_{\textrm{DOF}}}$ be an orthonormal basis for $U_h^K(\mathcal{T}_h)$. Then, $f_h(\bx,\bv_j)$ and $\rho_h(\bx)$ can be expanded as
\begin{equation}
f_h(\bx,\bv_j) = \sum_{i=1}^{N_{\textrm{DOF}}}\alpha_i^f(\bv_j)\phi_i(\bx)\quad\text{and}\quad\rho_h(\bx) = \sum_{i=1}^{N_{\textrm{DOF}}}\alpha_i^\rho \phi_i(\bx), \;\text{with}\; \alpha_i^\rho = \sum_{j=1}^{N_{\bv}}\omega_j\alpha_i^f(\bv_j).
\end{equation}
Define $\bff_j=\left(\alpha_1^f(\bv_j),\dots,\alpha_{N_{\textrm{DOF}}}^f(\bv_j)\right)^T\in\mathbb{R}^{N_{\textrm{DOF}}}$, $\brho=\left(\alpha_1^\rho,\dots,\alpha_{N_{\textrm{DOF}}}^\rho\right)^T\in\mathbb{R}^{N_{\textrm{DOF}}}$. Then, the DG scheme \eqref{eq:DG} can be written in its matrix-vector form:
\begin{equation}
    (\BD_j+\BSigma_t) \bff_j = \BSigma_s\brho+\bG+\bg_j^{(\textrm{bc})}=\BSigma_s\brho+\widetilde{\bG}_j,\;\brho=\sum_{j=1}^{N_{\bv}}\omega_j\bff_j,\quad j=1,\dots,N_{\bv},\label{eq:dg_matrix_vec}
\end{equation}
where $\BD_j, \BSigma_s, \BSigma_t\in\mathbb{R}^{N_{\textrm{DOF}}\times N_{\textrm{DOF}}}$, $\bG, \bg_j^{(\textrm{bc})}\in\mathbb{R}^{\pzc{N_\textrm{DOF}}}$ and $\widetilde{\bG}_j=\bG+\bg_j^{(\textrm{bc})}$ are defined as:
\begin{subequations}
\begin{align}
&(\BD_{j})_{kl} = -\sum_{i=1}^{N_x}\int_{\mathcal{T}_i} (\bv_j\cdot\nabla\phi_k(\bx)) \phi_l(\bx) d\bx+ \sum_{\mathcal{E}\in\partial\mathcal{T}_h\setminus\partial\mathcal{T}_{h,j}^{(\textrm{ibc})}}\int_{\mathcal{E}} \widehat{\BH}\left(\bv_j, \phi_l,\bn(\bx)\right)\phi_k(\bx) d\bx,
\\
&(\BSigma_t)_{kl} = \sum_{i=1}^{N_x}\int_{\mathcal{T}_i} \sigma_t(\bx)\phi_k(\bx) \phi_l(\bx) d\bx,
\quad (\BSigma_s)_{kl} = \sum_{i=1}^{N_x}\int_{\mathcal{T}_i} \sigma_s(\bx)\phi_k(\bx) \phi_l(\bx) d\bx,
\\
&(\bG)_k= \sum_{i=1}^{N_x}\int_{\mathcal{T}_i} G(\bx)\phi_k(\bx) d\bx,\quad
(\bg_j^{(\textrm{bc})})_k = -\sum_{\mathcal{E}\in\partial\mathcal{T}_{h,j}^{(\textrm{ibc})}}\int_{\mathcal{E}} g(\bx)\phi_k(\bx) \bv_j\cdot \bn(\bx) d\bx.
\end{align}
\end{subequations}
Let $\bff=(\bff_1^T,\dots,\bff_{N_{\bv}}^T)^T\in\mathbb{R}^{N_{\bv}N_{\textrm{DOF}}}$, then \eqref{eq:dg_matrix_vec} can be rewritten as a coupled linear system
\begin{equation}
\label{eq:one_equation}
\BA \bff =
\left(\begin{matrix}
\BD_1+\BSigma_t - \omega_1\BSigma_s & -\omega_2 \BSigma_s & \dots & -\omega_{N_{\bv}}\BSigma_s\\
-\omega_1\BSigma_s & \BD_2+\BSigma_t-\omega_2\BSigma_s & \dots & -\omega_{N_{\bv}}\BSigma_s \\
\vdots & \vdots & \vdots & \vdots \\
-\omega_1\BSigma_s & -\omega_2\BSigma_s & \dots & \BD_{N_{\bv}}+\BSigma_t-\omega_{N_{\bv}}\BSigma_s
\end{matrix}\right)
\left(
\begin{matrix}
\bff_1\\
\bff_2\\
\vdots\\
\bff_{N_{\bv}}
\end{matrix}
\right)
=
\left(\begin{matrix}
\wbG_1\\
\wbG_2\\
\vdots\\
\wbG_{N_{\bv}}
\end{matrix}
\right)=\bb
\end{equation}
\subsection{Source Iteration with Synthetic Acceleration}
We briefly review basic ideas of Source Iteration with Synthetic Acceleration (SISA) \cite{Adams2002FastIM}. Before delving into details, we present the  overall workflow of SISA in Alg. \ref{alg:SISA}. 
\begin{algorithm}[H]
\caption{Source Iteration with Synthetic Acceleration (SISA) to solve \eqref{eq:dg_matrix_vec}\label{alg:SISA} }
\label{alg:source}
\begin{algorithmic}
\STATE{Given initial guess: $\brho^{(0)}$, tolerance $\epssisa$, and the maximum number of iterations allowed $N_{\textrm{iter}}$.}
\STATE{\textbf{Source iteration}: solve $(\BD_{j}+\BSigma_t) \bff_{j}^{(1)} = \wbG_j + \BSigma_s \brho^{(0)}$ \pzc{with transport sweep}, 
$j=1,\dots,N_{\bv}$.}
\STATE{Update the density as $\brho^{(1,*)}=\sum_{j=1}^{N_{\bv}}\omega_j\bff_j^{(1)}$.}
\IF{ $||\rho^{(1,*)}-\rho^{(0)}||_\infty<\epsilon_{\textrm{SISA}}$,}
\STATE{return $\brho^{(1,*)}$ and $\bff^{(1)}$ as our solutions.}
\ENDIF
\STATE{\textbf{Correction}: solve a correction equation to obtain the correction $\delta\brho^{(1)}$.}
\STATE{Update the density as $\brho^{(1)}=\brho^{(1,*)}+\delta\brho^{(1)}$.}
\FOR{$l=2:N_{\textrm{iter}}$}
\STATE{\textbf{Source iteration}: solve $(\BD_j+\BSigma_t) \bff_{j}^{(l)} = \wbG_j + \BSigma_s \brho^{(l-1)}$ \pzc{with transport sweep}, $j=1,\dots,N_{\bv}$.}
\STATE{Update the density as $\brho^{(l,*)}=\sum_{j=1}^{N_{\bv}}\omega_j\bff_j^{(l)}$. }
\IF{ $||\rho^{(l,*)}-\rho^{(l-1)}||_\infty<\epsilon_{\textrm{SISA}}$,}
\STATE{return $\brho^{(l,*)}$ and $\bff^{(l)}$ as the solution.}
\ENDIF
\STATE{\textbf{Correction}: {solve a correction equation to obtain the correction $\delta\brho^{(l)}$.}}
\STATE{Update the density as $\brho^{(l)}=\brho^{(l,*)}+\delta\brho^{(l)}$.}
\ENDFOR
\end{algorithmic}
\end{algorithm}

\textbf{Source Iteration (SI):} In the linear system $\BA\in\mathbb{R}^{N_{\bv}N_{\textrm{DOF}}\times N_{\bv}N_{\textrm{DOF}}}$ (defiend in \eqref{eq:one_equation}), $\bff_j$'s are coupled through the integral term $\BSigma_s\brho=\sum_{j=1}^{N_{\bv}}\omega_j\BSigma_s \bff_j$. To avoid directly solving the huge coupled system $\BA$, SI decouples $\bff_j$ in each iteration by freezing the density $\brho$. Specifically, in the $l$-th source iteration $(l\geq 1)$, given the density determined by the previous iteration or the initial guess $\brho^{(l-1)}$, we update the solution by solving
\begin{align}
    (\BD_j+\BSigma_t)\bff_j^{(l)} = \BSigma_s\brho^{(l-1)}+\bG+\bg_j^{(\textrm{bc})},\; j=1,\dots,N_{\bv}.\label{eq:transport_sweep}
\end{align}
Without Synthetic Acceleration (SA), the density is updated as $\brho^{(l)}=\brho^{(l,*)}=\sum_{j=1}^{N_{\bv}}\omega_j\bff_j^{(l)}$.
With upwind DG spatial discretization, $\BD_j+\BSigma_t$ is a block lower triangular matrix, when a proper ordering of elements determined by the upwind direction for $\bv_j$ is applied. Consequently, the decoupled linear system \eqref{eq:transport_sweep} \pzc{is} efficiently solved in a matrix-free manner by a transport sweep \cite{Adams2002FastIM}.

\textbf{Synthetic Acceleration (SA):} It is well known that vanilla SI can suffer from arbitrarily slow convergence  \cite{Adams2002FastIM}. SA can be viewed as a preconditioning step to accelerate its convergence. SA introduces a correction to the density after each source iteration:
\begin{equation}
\brho^{(l)}=\brho^{(l,*)}+\delta\brho^{(l)},\quad\brho^{(l,*)}=\sum_{j=1}^{N_{\bv}}\omega_j\bff_j^{(l)}.
\end{equation}
Let $\bff_j$ be the exact solution to \eqref{eq:dg_matrix_vec}. The ideal density correction is 
\begin{equation}
\delta\brho^{(l)}=\sum_{j=1}^{N_{\bv}} \omega_j \delta\bff^{(l)},\quad \delta \bff_j^{(l)}=\bff_j-\bff_j^{(l)}.
\end{equation}
One can show that the ideal correction $\delta \bff_j^{(l)}$ satisfies the following discrete kinetic correction equation:
\begin{subequations}
\label{eq:discrete_correction}
\begin{align}
    &(\BD_{j}+\BSigma_t) \delta\bff_{j}^{(l)} = \BSigma_s \delta\brho^{(l)}+ \BSigma_s\left(\brho^{(l,*)}-\brho^{(l-1)}\right),\;
    j=1,\dots,N_{\bv},\\
    &\delta\brho^{(l)}=\sum_{j=1}^{N_{\bv}}\omega_j \delta\bff_j^{(l)},
\end{align}
\end{subequations}
which can be rewritten as 
\begin{equation}
\BA\delta\bff^{(l)} = \delta \bb^{(l)},\label{eq:one_equation_correction}
\end{equation}
where $\BA$ is defined in \eqref{eq:one_equation}, $\delta\bff^{(l)}=\left((\delta\bff_1^{(l)})^T,\dots,(\delta\bff^{(l)}_{N_{\bv}}
)^T\right)^T$, and 
$$\delta\bb^{(l)}=\underbrace{\left((\BSigma_s(\brho^{(l,*})-\brho^{(l-1)}))^T,\dots,(\BSigma_s(\brho^{(l,*)}-\brho^{(l-1)}))^T\right)^T}_{\text{Repeat } N_v\text{\;times}}.$$
Equation \eqref{eq:discrete_correction}
is the $S_N$-DG discretization of the kinetic correction equation 
\begin{subequations}
\label{eq:continuous_correction}
\begin{align}
    &\bv\cdot\nabla_{\bx}\delta f^{(l)}(\bx,\bv) +\sigma_t(\bx) \delta f^{(l)}(\bx,\bv) = \sigma_s\delta\rho^{(l)}(\bx)+\sigma_s (\bx)(\rho^{(l,*)}(\bx)-\rho^{(l-1)}(\bx)),\quad\bx\in\Omega_{\bx},\\
    &\delta\rho^{(l)}(\bx)=\frac{1}{4\pi}\int_{\mathbb{S}^2}\delta f^{(l)}(\bx,\bv)d\bv,\quad\bx\in\Omega_{\bx}\\
    &\delta f^{(l)}(\bx,\bv) = 0, \quad \bx\in \partial\Omega_{\bx},\;\bv\cdot \bn(\bx)<0,
\end{align}
\end{subequations}
which is a RTE with the isotropic source $\sigma_s(\bx)(\rho^{(l,*)}(\bx)-\rho^{(l-1)}(\bx))$ and zero inflow boundary conditions.
If one exactly solves the kinetic correction equation \eqref{eq:discrete_correction}, the SI converges with at most two iterations, because the solution of \eqref{eq:dg_matrix_vec} satisfies
\begin{equation*}
\brho=\brho^{(1)}=\brho^{(1,*)}+\delta\brho^{(1)}\quad\text{and}\quad\bff_j=\bff_j^{(2)}=(\BD_j+\BSigma_t)^{-1}(\BSigma_s\brho+\bG+\bg_j^{\textrm{(bc)}}).
\end{equation*} However, solving the kinetic correction equation \eqref{eq:discrete_correction} is as expensive as directly solving \eqref{eq:dg_matrix_vec}.

In practice, instead of the kinetic correction equation \eqref{eq:discrete_correction}, a computationally cheap low order approximation to it is solved in the correction step. 
For example, Diffusion Synthetic Acceleration (DSA) \cite{kopp1963synthetic,alcouffe1977dittusion,adams1992diffusion,wareing1993new,Adams2002FastIM} adopts the diffusion limit of the kinetic correction equation \eqref{eq:continuous_correction}
\begin{equation}
    -\nabla_{\bx}\cdot\left(\frac{1}{\sigma_s}\BD_{\bv}(\nabla_{\bx}\delta\rho^{(l)})\right) = -\sigma_a\delta\rho^{(l)}+\sigma_s(\rho^{(l,*)}-\rho^{(l-1)}), \label{eq:diffusion_correction}
\end{equation}
where $\BD_{\bv}=\textrm{diag}\left(\frac{1}{4\pi}\int_{\bv\in\mathbb{S}^2}\bv_x^2d\bv,\frac{1}{4\pi}\int_{\bv\in\mathbb{S}^2}\bv_y^2d\bv,\frac{1}{4\pi}\int_{\bv\in\mathbb{S}^2}\bv_z^2d\bv\right).$
Quasi-Diffusion method \cite{gol1964quasi,anistratov1993nonlinear,olivier2023family} uses the variable Eddington factor, and S$2$SA \cite{lorence1989s} employs a low order discrete ordinates ($S_N$) approximation. 
The effectiveness of these methods relies on the accuracy of their underlying low order approximations. For example, SI with DSA (SI-DSA) may converge slowly, if the kinetic correction equation is far from its diffusion limit \cite{Adams2002FastIM}. 
\subsubsection{\pzc{Memory efficient Source Iteration and Krylov subspace acceleration}\label{sec:krylov}}
\pzc{Utilizing  the fact that $\brho=\sum_{j=1}^{N_{\bv}}\omega_j\bff_j$, a more memory efficient Source Iteration and Krylov subspace solvers are developed. Here, we briefly review the basic ideas.}

\pzc{Solving equation \eqref{eq:dg_matrix_vec}, we obtain
\begin{equation}
    \bff_j = (\BD_j+\BSigma_t)^{-1} \BSigma_s\brho+(\BD_j+\BSigma_t)^{-1}(\bG+\bg_j^{(\textrm{bc})}),\quad\brho=\sum_{j=1}^{N_{\bv}}\omega_j\bff_j.\label{eq:f_matrix_vec}
\end{equation}
Numerically integrating $\bff$, equation \eqref{eq:f_matrix_vec} becomes
\begin{align}
\brho = \sum_{j=1}^{N_{\bv}}\left(\omega_j(\BD_j+\BSigma_t)^{-1}\right)\BSigma_s\brho + \sum_{j=1}^{N_{\bv}}\left(\omega_j(\BD_j+\BSigma_t)^{-1}\right)(\bG+\bg_j^{(\textrm{bc})}),\\
\left(\mathbf{I}-\sum_{j=1}^{N_{\bv}}\left(\omega_j(\BD_j+\BSigma_t)^{-1}\right)\BSigma_s\right)\brho = \sum_{j=1}^{N_{\bv}}\left(\omega_j(\BD_j+\BSigma_t)^{-1}\right)(\bG+\bg_j^{(\textrm{bc})}).
\label{eq:rho_matrix_vec}
\end{align}
Define 
\begin{equation}
\mathbf{L}=\sum_{j=1}^{N_{\bv}}\left(\omega_j(\BD_j+\BSigma_t)^{-1}\right)\BSigma_s\quad\text{and}\quad\overline{\bb}=\sum_{j=1}^{N_{\bv}}\left(\omega_j(\BD_j+\BSigma_t)^{-1}\right)(\bG+\bg_j^{\textrm{(bc)}}).
\end{equation} 
Equation \eqref{eq:rho_matrix_vec} can be rewritten as a linear system for $\brho$: 
\begin{equation}
    (\mathbf{I}-\mathbf{L})\brho=\overline{\bb}.
    \label{eq:rho_matrix_vec_def}
\end{equation} 
The matrix-vector product determined by $\mathbf{I}-\mathbf{L}=\mathbf{I}-\sum_{j=1}^{N_{\bv}}\left(\omega_j(\BD_j+\BSigma_t)^{-1}\right)\BSigma_s$ and $\mathbf{L}$ can be implemented in a matrix-free manner through transport sweep.}

\textbf{\pzc{Memory efficient Source Iteration}:}
\pzc{Utilizing \eqref{eq:rho_matrix_vec}, SISA can be seen as:}
\begin{equation}
\pzc{\brho^{(l,*)} = \mathbf{L}\brho^{(l-1)}, \quad \brho^{(l)}=\brho^{(l,*)}+\delta\brho^{(l)},}\label{eq:rho_SI}
\end{equation}
\pzc{where the correction $\delta\brho^{(l)}$ is obtained through Synthetic Acceleration.
Compared with \eqref{eq:transport_sweep}, the above formulation \label{eq:rho_SI} requires less memory, as one does not have to explicitly save $\bff_j$ for $j=1,\dots,N_{\bv}$. Though $\mathbf{L}\in\mathbb{R}^{N_{\textrm{DOF}}\times N_{\textrm{DOF}}}$, 
the computational cost of applying $\mathbf{L}$ is still $O(N_{\bv}N_{\textrm{DOF}})$, 
because it is realized in a matrix free manner through a transport sweep, which essentially solves linear systems determined by $\BD_j+\BSigma_t$ for all angular directions
$\bv_j, 1\leq j\leq N_{\bv}$ by sweeping along the corresponding upwind direction.}

\pzc{In summary, compared with the SI for $\bff$ in \eqref{eq:transport_sweep}, the SI for $\brho$ in \label{eq:rho_SI} requires less memory but shares the same leading order of computational cost due to transport sweeps.}  

\textbf{\pzc{DSA preconditioned Krylov solver}:}
\pzc{To accelerate the convergence of SI, Krylov methods such as the Generalized Minimal Residual method (GMRES) \cite{saad1986gmres} have been applied to solve $(\mathbf{I}-\mathbf{L})\brho=\overline{\bb}$.
However, GMRES without preconditioning is slower than SI-DSA in many cases \cite{guthrie1999gmres}. DSA based preconditioner for GMRES have been developed in \cite{warsa2004krylov} as follows. Under the SI framework, DSA obtains the density correction by   
\begin{equation}
\delta\brho^{(l)}=\mathbf{C}\BSigma_s(\brho^{(l,*)}-\brho^{(l-1)}),
\end{equation}
where operator $\mathbf{C}$ solves the linear system determined by the underlying discrete diffusion operator. Viewing SISA as a fixed point iteration for a left preconditioned linear system, \cite{warsa2004krylov} derives the DSA preconditioned linear system}
\begin{equation}
\pzc{(\mathbf{I}+\mathbf{C})(\mathbf{I}-\mathbf{L})\brho=(\mathbf{I}+\mathbf{C})\overline{\bb}.}
\end{equation}

\pzc{However, as demonstrated in \cite{ren2019fast}, even under GMRES framework, the DSA preconditioner may still become less effective without a sufficiently strong scattering effect \cite{ren2019fast}. Additionally, the DSA preconditioner does not utilize low rank structures with respect to parameters in parametric problems.}

\section{Reduced order model enhanced SISA\label{sec:romig-romsa}}
In applications such as uncertainty quantification, shape optimization and inverse problems, RTE needs to be solved many times for various parameters. 
Low rank structures with respect to the parameters of these parametric problems are not exploited by classical SA strategies like DSA. We propose utilizing data-driven reduced order models (ROMs) for these parametric problems and  their corresponding kinetic correction equations to enhance SISA by identifying and leveraging these low rank structures. 

The ROM for the parametric problem can be exploited to provide an improved initial guess for SISA. Additionally, the ROM for the kinetic correction equation can be exploited to design a low rank approximation to the kinetic correction equation \eqref{eq:discrete_correction}. Unlike empirical low order approximations such as the diffusion limit, this ROM-based low rank approximation directly builds on the kinetic description of the correction equation and leverages low rank structures with respect to the parameters. Furthermore, novel SA strategies can be developed based on this approximation.

This section is organized as follows. In Sec. \ref{sec:romig} we first build a ROM for parametric RTE following \cite{behne2023parametric} and then use the constructed ROM to provide an \pzc{improved} initial guess for SISA \pzc{or GMRES}. In Sec. \ref{sec:romsa-whole}, we propose a strategy to build a ROM for the kinetic correction equation without directly solving it, and then design two SA strategies built on this ROM for the kinetic correction equation.

\subsection{ROM based initial guess\label{sec:romig}}
We first construct a ROM for parametric RTE following \cite{behne2023parametric}, and then utilize the constructed ROM to provide an improved initial guess for the SISA.

\subsubsection{ROM for parametric RTE\label{sec:rom-rte}}
Following \cite{behne2023parametric}, we employ Proper Orthogonal Decomposition (POD) \cite{sirovich1987turbulence,berkooz1993proper} to build a ROM for parametric RTE. POD applies an offline-online decomposition framework. In the online stage, we seek reduced order solutions in a reduced order space constructed offline. The dimension of this reduced order space is much smaller than $N_vN_{\textrm{DOF}}$, resulting in very efficient prediction online. In the offline stage, we construct the reduced order space by extracting low rank structures from the training data. The details of this procedure are as follows.

\textbf{Offline reduced order space construction:} 
Let $\mathcal{P}_{\textrm{train}}=\{\bmu_j\in\mathbb{R}^p\}_{j=1}^{N_{\textrm{train}}}$ be a training set of parameters, such as the strength of scattering and absorption cross sections, boundary conditions and  source terms.
Denote the solution corresponding to a parameter $\bmu$ as $\bff^{\bmu}$, also called a snapshot. The snapshot matrix determined by the training set is defined as
\begin{equation}
\BF=\left(\bff^{\bmu_1},\bff^{\bmu_2},\dots,\bff^{\bmu_{N_{\textrm{train}}}}\right).
\end{equation}
Each column of the snapshot matrix $\BF$ is a full order solution corresponding to a parameter $\bmu\in\mathcal{P}_{\textrm{train}}$. These solutions can be obtained by solving \eqref{eq:dg_matrix_vec} with SI-DSA.

To identify low rank structures in the training data, we compute the singular value decomposition (SVD) of the snapshot matrix $\BF$:
\begin{equation}
    \BF = \BU \BLambda \BV^T,
\end{equation}
where $\BU\in\mathbb{R}^{N_{\bv}N_{\textrm{DOF}}\times N_{\bv}N_{\textrm{DOF}}}$, $\BV\in\mathbb{R}^{N_{\textrm{train}}\times{N_{\textrm{train}}} }$ are orthogonal matrices, and the $k$-th diagonal element of $\BLambda\in\mathbb{R}^{N_{\bv}N_{\textrm{DOF}}\times N_{\textrm{train}}}$, $\lambda_k$, is the $k$-th largest singular value of $\BF$.
Given an error tolerance $\epsilon_{\textrm{POD}}$, we determine the dimension of the desired reduced order space $r_p$ as
\begin{equation}
    r_p = \arg\min_{l} \frac{\sum_{k=1}^l{\lambda_k}}{\textrm{trace}(\BLambda)}\geq 1-\epsilon_{\textrm{POD}}. \label{eq:pod_truncation}
\end{equation}
The POD basis $\BU_{r_p}$ is defined as the first $r_p$ columns of $\BU$. The constructed reduced order space is the space determined by columns of $\BU_{r_p}$. In the online stage, we seek reduced order solutions in this space. 

\textbf{Online prediction:} In the online stage, to efficiently predict the solution for  a parameter $\bmu\not\in\mathcal{P}_{\textrm{train}}$, we seek a reduced order solution in the reduced order space, $\bff^{\bmu}\approx \BU_{r_p}\bc_{\bmu,r_p}$. To determine the coefficient of the reduced order solution $\bc_{\bmu,r_p}$, we project the full order problem \eqref{eq:one_equation} onto the reduced order space through a Galerkin projection:
\begin{equation}
\BA_{\bmu,r_p}\bc_{\bmu,r_p} = \BU_{r_p}^T \BA_{\bmu}\BU_{r_p} \bc_{\bmu,r_p} = \BU_{r_p}^T \bb_{\bmu}.
\label{eq:rom_parametric}
\end{equation}
Here, footnote $\bmu$ represents $\bmu$-dependence. 

\textbf{Offline precomputations: }
We can utilize offline precomputations to efficiently construct the  reduced order operator $\BA_{\bmu,r_p}$ online. 
For now, we consider affine operators $\BA_{\bmu}$ satisfying
\begin{equation}
\BA_{\bmu} = \sum_{k=1}^{m}\psi_k(\bmu)\BA_k,\; \psi_k(\bmu)\in\mathbb{R},\;\BA_k\in\mathbb{R}^{N_{\bv}N_{\textrm{DOF}}\times{N_{\bv}N_{\textrm{DOF}}}},
\end{equation}
where $\BA_k$'s are constant matrices independent of the parameter $\bmu$. Under this affine assumption, we can precompute and save reduced operators $\BA_{k,r_p}=\BU_{r_p}^T\BA_k \BU_{r_p}\in\mathbb{R}^{{r_p}\times {r_p}}$ offline. In the online stage, the reduced order operator $\BA_{\bmu,r_p}$ can be efficiently constructed for a new parameter $\bmu\not\in\mathcal{P}_{\textrm{train}}$ as $\BA_{\bmu,r_p}=\sum_{k=1}^m \psi_k(\bmu) \BA_{k,r_p}$. 
When $\BA_{\bmu}$ has non-affine dependence on the parameter $\bmu$, one needs to apply empirical interpolation method (EIM) or discrete EIM (DEIM) \cite{barrault2004empirical,chaturantabut2010nonlinear} to achieve significant computational saving in the online stage. The main focus of this paper is to demonstrate the potential of using ROMs to accelerate SISA. Applying EIM or DEIM for non-affine problems will be left for future investigations. 

Furthermore, offline precomputations can help us to compute the density for the reduced order solution or project isotropic source terms more efficiently online. 
Let $\BU_{r_p,j}$ be the submatrix of $\BU_{r_p}$ corresponding to its $(j-1)N_{\textrm{DOF}}+1$-th row to its $jN_{\textrm{DOF}}$-th row, i.e., rows aligned with $\bff_{j}$. Then, $\BU_{r_p} = (\BU^T_{r_p,1},\dots,\BU^T_{r_p,N_{\bv}})^T$. 
The reduced order approximation to the density is given by 
\begin{equation}
    \brho^\bmu\approx \sum_{j=1}^{N_{\bv}} \omega_j\BU_{r_p,j}\bc_{\bmu,r_p} =  \BU_{r_p}^{\rho}\bc_{\mu,r_p},\quad \BU_{r_p}^\rho = \sum_{j=1}^{N_{\bv}}\omega_j\BU_{r_p,j}.
\end{equation} 
When the source term in RTE is isotropic and the boundary conditions are zero inflow boundary conditions, we have $\bb_{\bmu}=(\bG^T_{\bmu},\dots,\bG^T_{\bmu})^T$ and
\begin{equation}
\BU_{r_p}^T\bb_\bmu = \sum_{j=1}^{N_{\bv}} \BU_{r_p,j}^T\bG_{\bmu} = \BU_{r_p}^{\textrm{iso}} \bG_{\bmu}, \; \BU_{r_p}^{\textrm{iso}}=\sum_{j=1}^{N_{\bv}}U_{r_p,j}.
\end{equation}
To improve online efficiency, we precompute and save $\BU^\rho_{r_p}$ and $\BU_{r_p}^{\textrm{iso}}$ offline:
\begin{equation}
    \BU_{r_p}^\rho = \sum_{j=1}^{N_{\bv}}\omega_j\BU_{r_p,j}\quad\text{and}\quad
    \BU_{r_p}^{\textrm{iso}} = \sum_{j=1}^{N_{\bv}} \BU_{r_p,j}.
    \label{eq:rho_basis}
\end{equation}
.

\textbf{ROM based initial guess (ROMIG):} The ROM for the parametric problem \eqref{eq:rom_parametric} can be seen as a surrogate solver to predict solutions for parameters outside the training set. The size of the reduced order problem \eqref{eq:rom_parametric} is $r_p$, where $r_p\ll N_{\textrm{DOF}}N_{\bv}$. Consequently, solving this reduced order problem  is highly efficient. Furthermore, the reduced order solution can be employed as an initial guess for SISA \pzc{or Krylov solvers in Sec. \ref{sec:krylov}}. This ROM-based initial guess (ROMIG) is usually closer to the exact solution compared with a random or $0$ initial guess. 
We want to mention that, beyond the scope of RTE, a ROM based initial guess has been applied to \pzc{enhance Conjugate Gradient solvers with multigrid preconditioners for parametric} elliptic equations \cite{hou2023reduced}.

\begin{rem}
\pzc{The current offline strategy can be further enhanced. First, in our numerical tests, we uniformly sample parameters in a training set to construct reduced basis and reduced operators. This sampling strategy may lead to unnecessary computational cost due to over-sampling or errors due to under-sampling. To address this issue, we can apply a greedy algorithm \cite{hesthaven2016certified} based on error estimators to adaptively sample more informative training parameters. Such a greedy strategy allows adaptively sampling more parameters in parameter regions with sharper features in the solution manifold and fewer parameters in smoother regions.
Furthermore, incremental SVD \cite{brand2002incremental} or randomized SVD \cite{martinsson2020randomized} can be leveraged to more efficiently generate reduced order basis.}    
\end{rem}

\begin{rem}
\pzc{When parameters vary dramatically in their magnitudes, the dimension of the reduced order space may become relatively high for such cases. To balance online efficiency and accuracy, an adaptive domain decomposition in the parameter space can be applied to construct localized ROMs in each subdomain of the parameter space \cite{eftang2010hp,maday2013locally}.}  
\end{rem}

\begin{rem}\label{rem:minimal-intrusive}
If only the source term and the boundary conditions rely on the parameter $\bmu$,  a more efficient ROM requiring less memory can be built for parametric RTE \cite{behne2022minimally} \pzc{based on the memory efficient formulation in Sec. \ref{sec:krylov}}. However, when the scattering or absorption cross section depends on the parameter $\bmu$, whether the dependence is affine or not, the online stage of the ROM in \cite{behne2022minimally} may become \pzc{significantly less efficient}. 
\end{rem}
\subsection{ROM based Synthetic Acceleration\label{sec:romsa-whole}}
The key step of SA is to correct the density by solving a computationally cheap low order approximation to the kinetic correction equation \eqref{eq:discrete_correction}. Classical SA strategies apply empirical low order approximations such as the diffusion limit. They become less effective when their underlying low order approximations are not accurate enough. Additionally, they do not exploit low rank structures with respect to the parameters of parametric problems.

Alternatively, we propose to adopt a low rank approximation to the kinetic correction equation based on a data-driven ROM for it. This low rank approximation directly builds on the kinetic description of the correction equation and  leverages low rank structures with respect to the parameters.

To construct a data-driven ROM for the kinetic correction equation \eqref{eq:discrete_correction}, solutions to it are required. However, directly solving the kinetic correction equation 
for various parameters $\bmu$ and source terms for multiple iterations $\BSigma_{s,\bmu}(\brho^{\bmu,(l,*)}-\brho^{\bmu,(l-1)})$ ($l\geq1$) dramatically increases the computational cost of the offline stage.  Fortunately, when solving RTE with SISA, solutions to the kinetic correction equations for source terms corresponding to multiple iterations can be obtained without directly solving them. The approach of accomplishing this will be detailed in Sec. \ref{sec:rom-correction}. Subsequently, in Sec. \ref{sec:romsa}, we propose two SA strategies based on the constructed ROM for the kinetic correction equation.

\subsubsection{ROM for the correction equation\label{sec:rom-correction}}
Given the snapshot matrix for the kinetic correction equation\eqref{eq:discrete_correction}, $\delta\BF$, we can build a ROM for it following a procedure similar to Sec. \ref{sec:romig}. We denote the resulting POD basis as $\delta\BU_{r_c}\in\mathbb{R}^{N_{\textrm{DOF}}\times r_c}$, and the dimension of the corresponding reduced order space as $r_c$.
Below, we outline an efficient approach to generate the snapshot matrix without directly solving the kinetic correction equation. 

\textbf{Generation of the snapshot matrix:} We choose a window size $\mfw\in\mathbb{Z},\;\mfw\geq 1$. When solving parametric RTE with SISA, we save not only the converged solution $\bff^{\bmu}$ but also the intermediate solutions at the $l$-th iteration $\bff^{\bmu,(l)}$ for $1\leq l\leq \mfw_{\bmu}$, where 
\[
\mfw_{\bmu}=\min(n_{\bmu}^{(\textrm{conv})},\mfw)
\]
and $n_{\bmu}^{(\textrm{conv})}$ is the number of iterations for convergence. Recall that the solution of the kinetic correction equation \eqref{eq:discrete_correction} for the $l$-th iteration is $\delta \bff^{(l)}=\bff-\bff^{(l)}$.
Consequently, the snapshot matrix for the kinetic correction equation $\delta\BF$
can be assembled as 
\begin{equation}
\delta\BF = \left(\delta \bff^{\bmu_1,(1)},\dots,\delta\bff^{\bmu_1,(\mfw_{\bmu_1})},\delta \bff^{\bmu_2,(1)},\dots,\delta\bff^{\bmu_2,(\mfw_{\bmu_2})},\delta\bff^{\bmu_{N_{\textrm{train}}},(1)} ,\dots,\delta\bff^{\bmu_{N_{\textrm{train}}},(\mfw_{\bmu_{N_{\textrm{train}}}})} \right),
\end{equation}
where $\delta\bff^{\bmu_k,(l)}=\bff^{\bmu_k}-\bff^{\bmu_k,(l)}$. In summary, to construct the snapshot matrix $\delta\BF$, it is sufficient to save both the converged solution $\bff^{\bmu}$ and the intermediate results for the first few iterations $\bff^{\bmu,(l)}$, $l=1,\dots,\mfw_{\bmu}$.

\textbf{Online computation:} Following Sec. \ref{sec:rom-rte}, we can build the POD basis $\delta\BU_{r_c}$ based on the SVD of the snapshot matrix $\delta\BF$.
For a parameter $\bmu\not\in\mathcal{P}_{\textrm{train}}$, we project the kinetic correction equation for the $l$-th iteration \eqref{eq:one_equation_correction} to the $r_c$ dimensional reduced order space determined by $\delta \BU_{r_c}$, and seek a reduced order solution  $\delta \bff^{\bmu,(l)}\approx \delta \BU_{r_c} \delta \bc^{(l)}_{\bmu,r_c}$ by solving
\begin{equation}
   \delta \BU_{r_c}^T \BA_{\mu}\delta \BU_{r_c} \delta \bc^{(l)}_{\bmu,r_c} = \delta \BU_{r_c}^T\delta\bb^{(l)}_{\bmu}.
   \label{eq:rom_correction}
\end{equation}
Similar to \eqref{eq:rho_basis}, we can define $\delta \BU_{r_c}^\rho$ and $\delta \BU_{r_c}^{\textrm{iso}}$. Then, the reduced order density correction is $\delta\brho^{\bmu,(l)}\approx \delta\BU_{r_c}^\rho\delta\bc^{(l)}_{\bmu,r_c}$. Moreover, the kinetic correction equation \eqref{eq:continuous_correction} always has zero inflow boundary conditions and isotropic source term, hence its right hand side can always be projected as 
$(\delta\BU_{r_c}^{\textrm{iso}})^T\left(\BSigma_{s,\mu}(\brho^{\bmu,(l,*)}-\brho^{\bmu,(l-1)})\right)$.

\begin{figure}[]
  \begin{center} 
\includegraphics[width=0.31\textwidth]{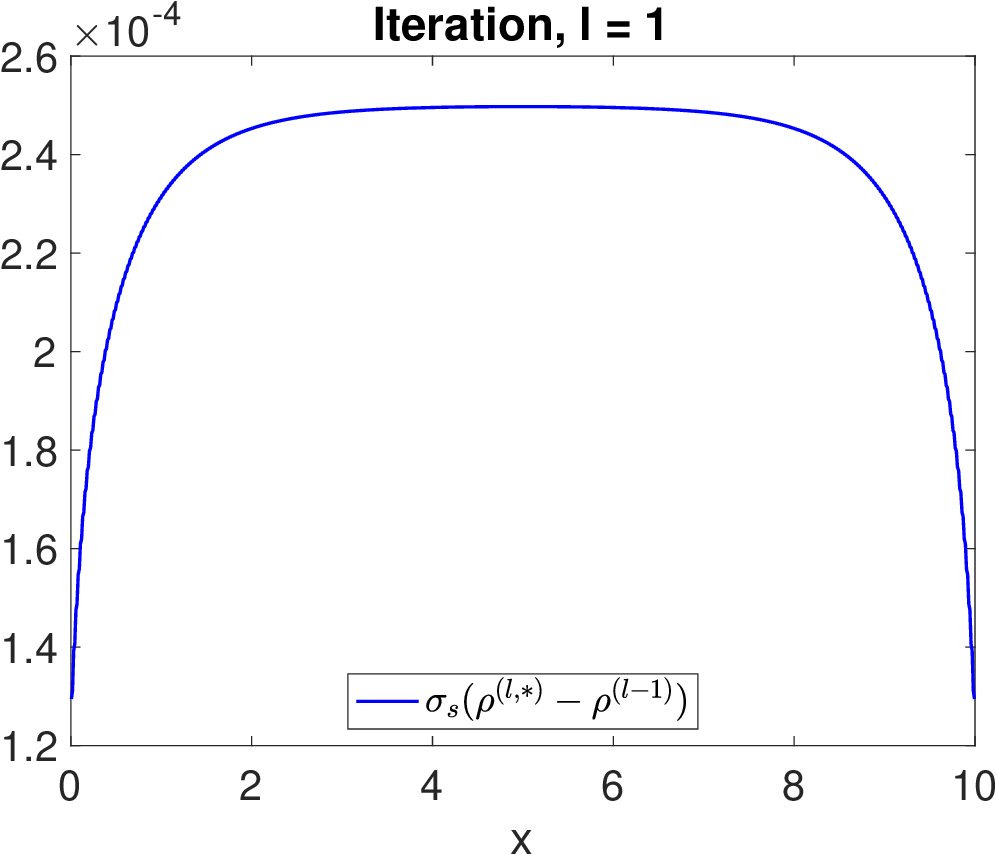}
\includegraphics[width=0.31\textwidth]{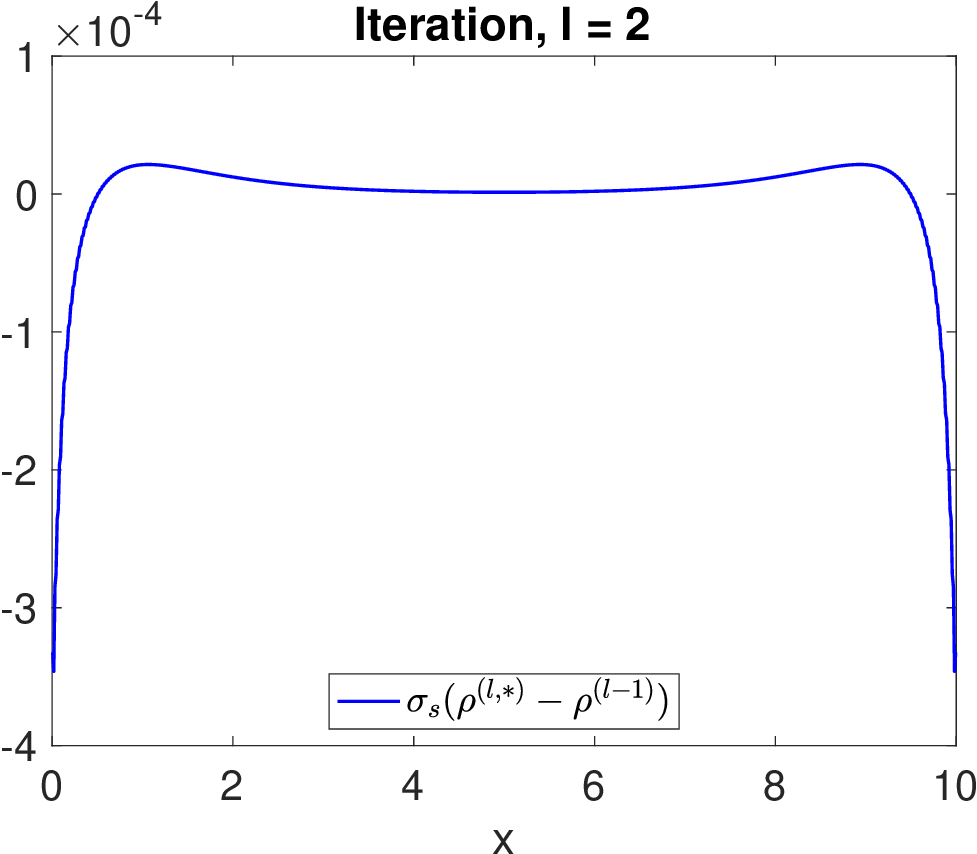}
\includegraphics[width=0.31\textwidth]{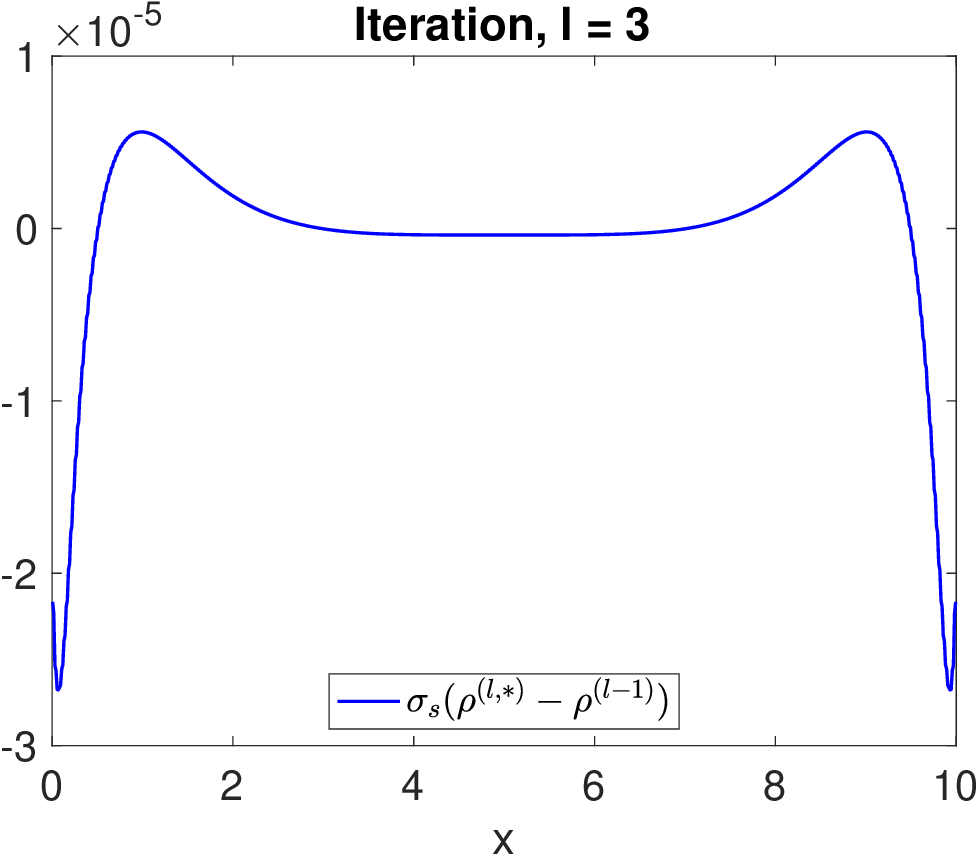}
\includegraphics[width=0.31\textwidth]{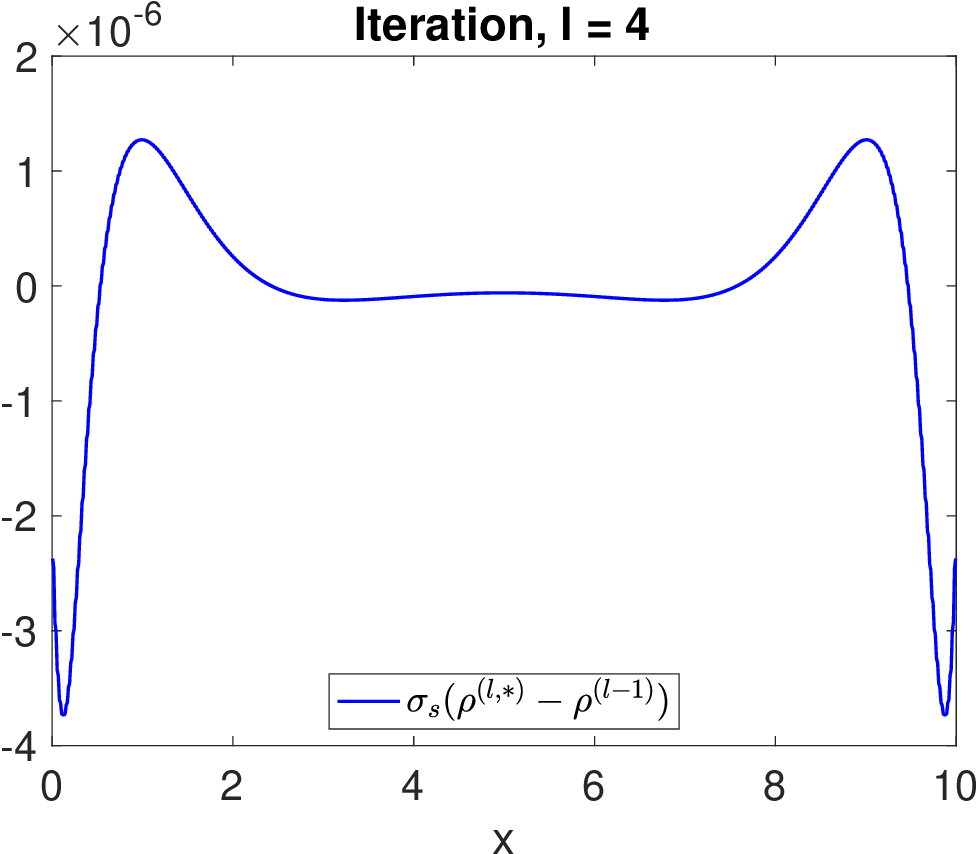}
\includegraphics[width=0.31\textwidth]{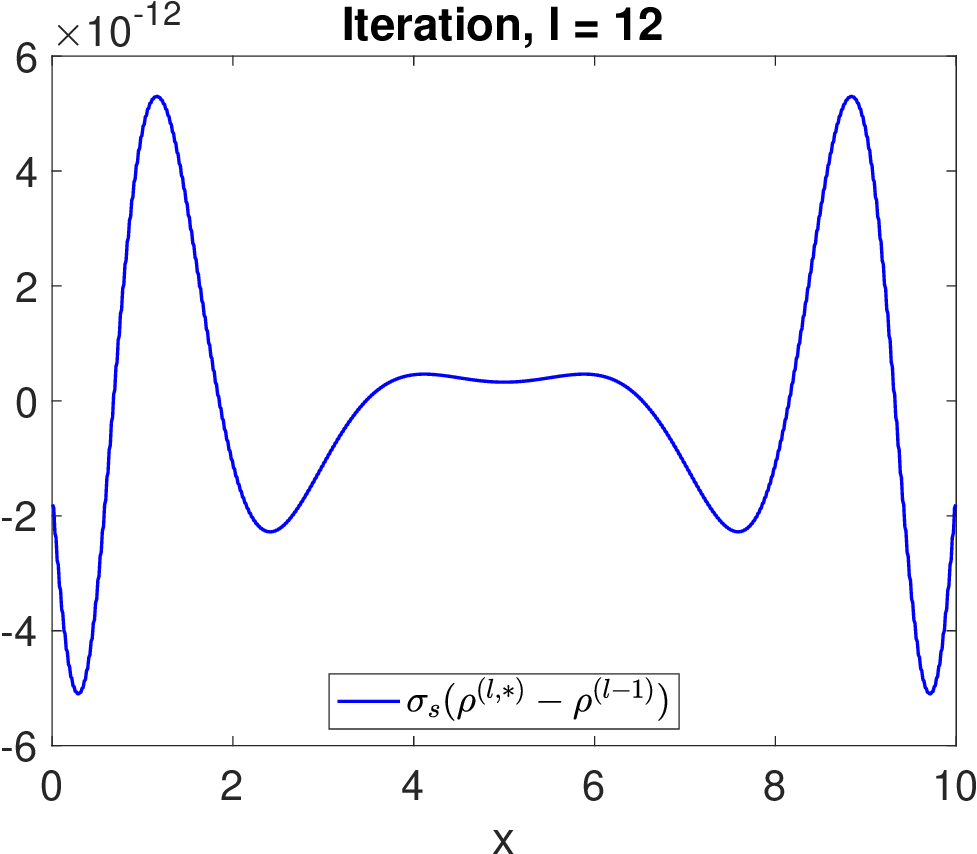}
\includegraphics[width=0.31\textwidth]{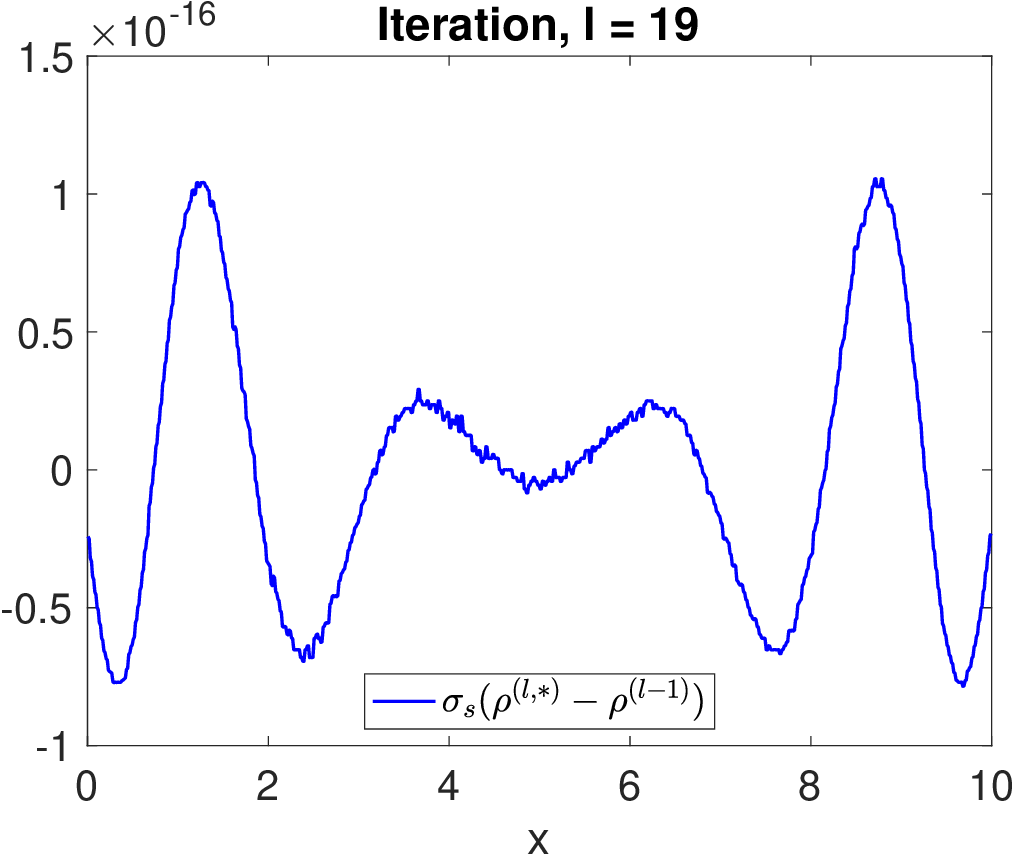}
  \caption{The source term in the correction equation for the $l$-th ($l=1,2,3,4,12,19$) iteration in the SI-DSA solving the problem in Sec. \ref{sec:rom-correction}. \label{fig:source_trajectory}}
  \end{center}
\end{figure}

\textbf{The role of the window size $\mfw$:}  As source iterations continue, the isotropic source in the correction equation for the $l$-th iteration, i.e. $\BSigma_s(\brho^{(l,*)}-\brho^{(l-1)})$, changes its shape. To demonstrate this shape variation, we consider the following example in 1D slab geometry: 
\begin{align*}
\Omega_{\bx}=[0,10] ,\;G = 0.01, \;
\sigma_t(x)=\sigma_s(x)=1,\;
f(0,\xi) = 0\;\text{with } \xi>0, \;f(10,\xi) = 0 \;\text{with } \xi<0.  
\end{align*}
We use $16$ Gauss-Legendre points for angular discretization, linear upwind DG for spatial discretization and a uniform mesh with $400$ elements to partition $\Omega_{\bx}$. The resulting linear system is solved with SI-DSA. In Fig. \ref{fig:source_trajectory}, we present the shape of the source in the correction equation for the $l$-th iteration, i.e. $\sigma_s(x)(\rho^{(l,*)}(x)-\rho^{(l-1)}(x))$, where $l=1,2,3,4,12,19$. We observe that the isotropic source term becomes more oscillatory as iterations continue. To account for this shape variation of the source term, we use a window size $\mfw\geq1$. 

The larger window size $\mfw$ is, the more information for later iterations will be included in our training data. However, in practice, we can only use a moderate window size $\mfw$ to avoid unacceptable memory costs. On the other hand, when the tolerance in the stopping criteria of SISA is small, many iterations may be needed and the shape of source terms may vary significantly. Nevertheless, using a moderate window size $\mfw$, our ROM may not be accurate enough in the later stage of SISA. We will discuss how to handle this issue when designing a ROM-based SA strategy in Sec. \ref{sec:romsa}.

\subsubsection{Two ROM based SA strategies\label{sec:romsa}}
We design a ROM-based SA strategy (ROMSA) and a hybridized strategy called ROMSAD, which combines ROMSA and DSA.

\textbf{ROMSA:} If the underlying ROM is generated with a window size $\mfw$, we call our SA strategy, ROMSA-$\mfw$. ROMSA uses a ROM-based low rank approximation to the kinetic correction equation \eqref{eq:discrete_correction}. Specifically,
after the $l$-th source iteration, ROMSA solves the reduced order kinetic correction equation \eqref{eq:rom_correction} to get the density correction $\delta\brho^{\bmu,(l)}= \delta\BU_{r_c}^\rho\delta\bc^{(l)}_{\bmu,r_c}$, where $\bc^{(l)}_{\bmu,r_c}$ is the solution of \eqref{eq:rom_correction}.
This ROM-based low rank approximation directly builds on the kinetic description of the correction equation and leverages low rank structures with respect to the parameters of the underlying parametric problem.

\begin{figure}[]
  \begin{center} 
\includegraphics[width=0.45\textwidth]{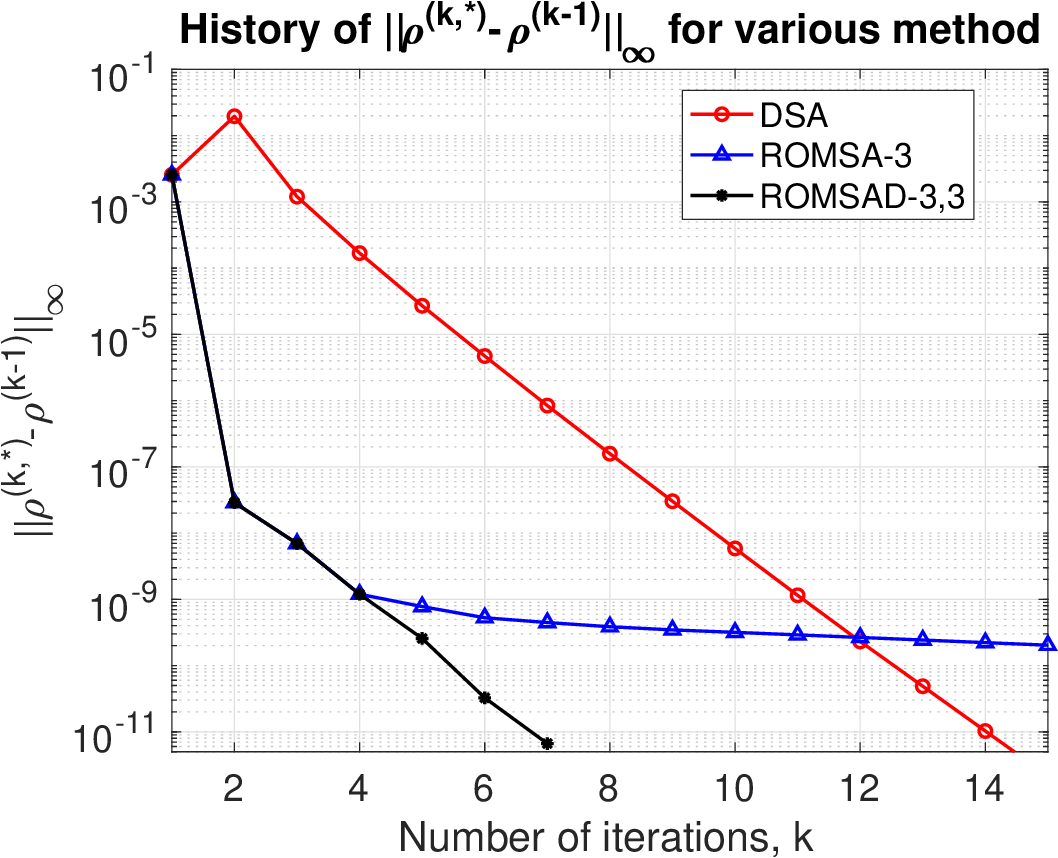}
  \caption{The history of the change of $\brho$ in the $k$-th iteration, i.e. $||\brho^{(k,*)}-\brho^{(k-1)}||_\infty$, for  the cross-regime problem in Sec. \ref{sec:num-1d-cross} with $\mu_s=19.772$ and $\epspod=10^{-9}$. \label{fig:convergence_demo}}
  \end{center}
\end{figure}

Numerically, we observe that ROMSA-$\mfw$ achieves greater acceleration than DSA in the early stage of SISA. However, as iterations continue, the shapes of the source terms in the correction equations may vary significantly. With a moderate window size $\mfw$, the underlying ROM may fail to accurately approximate the kinetic correction equation in the later stage of SISA. Consequently, ROMSA may suffer from an efficiency reduction in the later stage.  
For a clear illustration, in Fig. \ref{fig:convergence_demo}, we present how the difference of $\brho$ for successive iterations, $||\brho^{(k,*)}-\brho^{(k-1)}||_\infty$, evolves for the cross-regime problem in Sec. \ref{sec:num-1d-cross}. We observe a rapid decrease in $||\brho^{(k,*)}-\brho^{(k-1)}||_\infty$ for ROMSA-$3$ in the first $4$ iterations, followed by a slow decline starting from the $5$-th iteration. The decrease rate of $||\brho^{(k,*)}-\brho^{(k-1)}||_\infty$ for DSA is slower than or comparable to the rate for ROMSA in the first $4$ iterations, but becomes faster than  ROMSA from the $5$-th iteration and onward. 

In summary, with a moderate window size $\mfw$, ROMSA is more efficient than DSA in the first few iterations of SISA, but less robust than DSA. Hence, we should only use ROMSA in the early stage when it is more efficient, and then switch to DSA to avoid the potential efficiency reduction of ROMSA.

\begin{rem}
In our numerical tests, we observe that, for some parametric problems, SI accelerated by ROMSA converges very fast and does not suffer from the efficiency decline. Nevertheless, the lack of robustness severely restricts the applicability of ROMSA for general parametric problems. In contrary, DSA is robust for all tests. 
\end{rem}

\textbf{ROMSA hybridized with DSA (ROMSAD):} To leverage high efficiency of ROMSA in the early stage of SISA and the robustness of DSA, we propose a hybridized strategy called ROMSAD.

As observed in Fig. \ref{fig:convergence_demo}, ROMSA is highly efficient in the early stage of SISA, while DSA maintains robust across all stages. We propose to use ROMSA in the first few iterations, and then automatically switch to DSA to avoid the potential efficiency reduction of ROMSA. Specifically, we use ROMSA after the $l$-th source iteration if
$$||\rho^{(l,*)}-\rho^{(l-1)}||_{\infty}\geq \epsilon_{\textrm{ROMSAD}}\quad\text{and}\quad l\leq \theta.$$ Otherwise, we use DSA.
We choose 
\begin{equation}
 \epsromsad = \eta\max(\epssisa^{\textrm{train}},\epspod),\quad \eta\geq 0, \label{eq:romsad_parameters}   
\end{equation}
where $\epssisa^{\textrm{train}}$ is the tolerance in the stopping criteria of SISA generating our training data and $\epspod$ is the truncation tolerance in POD (see equation \eqref{eq:pod_truncation}). These two values determine the accuracy of the ROM utilized by ROMSA. 
We summarize our ROMSAD algorithm in Alg. \ref{alg:romsad}.
We refer to  ROMSAD as ROMSAD-$\mfw,\theta$, when the underlying ROM is generated with a window size $\mfw$, and ROMSA is allowed to be applied after at most  the first $\theta$ source iterations. We denote SI with ROMSAD as SI-ROMSAD.

To demonstrate the effectiveness of ROMSAD, in Fig. \ref{fig:convergence_demo}, we also show how $||\brho^{(k,*)}-\brho^{(k-1)}||_\infty$ evolves for ROMSAD-$3,3$. ROMSAD-$3,3$ corrects the density with ROMSA after the first $3$ source iterations, and then switches to DSA. As a result, ROMSAD enjoys a decrease of $||\brho^{(k,*)}-\brho^{(k-1)}||_\infty$ as rapid as ROMSA in the first $4$ iterations, and thereafter demonstrates performance as robust as DSA. For this problem, overall, ROMSAD converges faster than both ROMSA and DSA.  ROMSAD resolves the robustness issue of ROMSA without sacrificing its great efficiency in the early stage of SISA.

\begin{algorithm}[H]
\caption{Determine applying ROMSA or DSA after the $l$-th source iteration of SI-ROMSAD. \label{alg:romsad} }
\label{alg:source}
\begin{algorithmic}
\STATE{Input: current number of iterations $l$, a window size $\theta$, and a tolerance $\epsilon_{\textrm{ROMSAD}}$.}
\IF{$1\leq l\leq \theta$ and $||\rho^{(l,*)}-\rho^{(l-1)}||_{\infty}\geq \epsilon_{\textrm{ROMSAD}}$,}
    \STATE{apply ROMSA.}
\ELSE
    \STATE{apply DSA.}
\ENDIF
\end{algorithmic}
\end{algorithm}

\begin{rem}\label{rem:challenge_romigsa}
At first glance, it seems to be natural to directly use a ROM-based initial guess (ROMIG) for SI-ROMSAD. However, when using a ROMIG close to the exact solution, the shapes of source terms in the correction equations can differ significantly from those for other initial guesses, even in the early stage of SISA. Consequently, the underlying ROM of ROMSAD may not be accurate, unless it is built based on training data generated with the ROMIG.

A potential solution is to employ an iterative offline stage that adaptively updates ROMs for the parametric problem and the kinetic correction equation. Specifically, we first build an initial ROM for the parametric problem and then use an initial guess based on this ROM to generate training data for building an initial ROM for the correction equation. After initializing these ROMs, we can adaptively update them in an alternating manner. This iterative offline strategy will be left for future investigation.
\end{rem}

\section{Numerical results\label{sec:numerical}}
In our numerical tests, we consider both  1D slab geometry and 2D $X$-$Y$ geometry. 
\pzc{We compare the performance of the proposed  methods with both SI-DSA and DSA preconditioned GMRES solver. }
\pzc{Numerical tests in 1D slab geometry aim to investigate the influence of various parameters in the proposed method, while numerical tests in 2D X-Y geometry target on comparing the computational efficiency of various methods.}

In all examples, we use linear upwind DG spatial discretization $(K=1$ in \eqref{eq:discrete_space}).  We set $\eta$  as $\eta=0.1$ in the threshold of the switching criteria of ROMSAD (see \eqref{eq:romsad_parameters}).  We refer to the training set for parameters as $\mathcal{P}_{\textrm{train}}$, and the test set for parameters as $\mathcal{P}_{\textrm{test}}$. The training set and the test set satisfy $\mathcal{P}_{\textrm{train}}\cap\mathcal{P}_{\textrm{test}}=\emptyset$.  

\subsection{Implementation details and notations}
\pzc{Here, we summarize the essential implementation details of our code and the notations used throughout this section.}

\subsubsection{Implementation details \label{sec:details}}
\pzc{\textbf{DSA:}} We use a consistent \pzc{DSA} following \cite{larsen1982unconditionally,Adams2002FastIM}. \pzc{The discretization of the diffusion correction equation \eqref{eq:diffusion_correction} is derived by first substituting the $P_1$ ansatz $\delta f(\bx,\bv)=\delta\rho(\bx)+3\bv \delta J(\bx,\bv)$ into the upwind DG discretization of the kinetic correction equation \eqref{eq:continuous_correction} and then numerically taking the zero-th and first order moment of the resulting discrete system in the angular space (see Appx. \ref{sec:dsa} for more details).} In 1D slab geometry,  the resulting \pzc{discrete} diffusion equation is solved by a direct solver. In 2D $X$-$Y$ geometry, the \pzc{discrete} diffusion equation is solved by an algebraic multigrid (AMG) solver implemented based on the iFEM package \cite{Chen:2008ifem}. 

\pzc{\textbf{ROMSA:}} The \pzc{size of} reduced order kinetic correction problem for ROMSA is always \pzc{very small}, \pzc{so it is} solved by a direct solver.

\pzc{\textbf{Preconditioned GMRES:} As outlined in Sec. \ref{sec:krylov}, we follow the ideas in \cite{warsa2004krylov} to implement the DSA preconditioned GMRES solver with restart. Specifically, we provide function handles, which realize the operation of DSA preconditioner utilizing AMG and the matrix-vector product $(\mathbf{I}-\mathbf{L})\brho$ through transport sweeps, to Matlab's GMRES solver.}

\pzc{We want to point out that the number of transport sweeps required for convergence is slightly more than the number of iterations in GMRES. In Alg. \ref{alg:gmres} of Appx. \ref{sec:gmres}, we outline Matlab's GMRES implementation, which follows \cite{barrett1994templates}. The sources of extra transport sweeps are as follows. 
} 
\begin{enumerate}
    \item \pzc{One transport sweep is required to construct the right hand side of \eqref{eq:rho_matrix_vec}, which is defined in \eqref{eq:rho_matrix_vec_def}.}
    \item \pzc{One transport sweep is needed for the initialization (line 2 of Alg. \ref{alg:gmres}) and (line 12 of Alg. \ref{alg:gmres}) re-initialization before restart.} 
\end{enumerate}
\pzc{We set the number of iterations before restart as $25$, and in all of our numerical tests DSA-preconditioned GMRES converges before restart. As a result, the number of transport sweeps required by GMRES solver minus the number of iterations for its convergence is at most $3$ in our tests.}

\pzc{\textbf{Initial guess for iterative solvers:} unless using ROMIG, SISA and GMRES always starts with zero initial guess.}

\pzc{\textbf{Offline construction of ROMs and online test}:} \pzc{During the offline stage of constructing the ROM, we employ SI-DSA  to generate the snapshot matrix for $f$, directly solving the equation for $\bff_j$ \eqref{eq:transport_sweep} in each iteration.} 

\pzc{In the online test, we apply the memory efficient formulation of SI and Krylov solver outlined in Sec. \ref{sec:krylov}. It is important to note that the leading order computational cost of memory efficient SI for $\brho$ in \eqref{eq:rho_SI}  is the same as the original formulation for $\bff_j$ in \eqref{eq:transport_sweep}. In our numerical tests, we observe that the computational time of these two formulations of SI are nearly identical. }

\subsubsection{Notations \label{sec:numerical-notations}}
\pzc{Here, we introduce notations used throughout this section, and summarize them in Tab. \ref{tab:notation}. We will use ``DSA" to represent SI-DSA using zero initial guess, ``ROMIG" to represent SI-DSA using ROM based initial guess, ``ROMSAD" to represent SI-ROMSAD, ``PGMRES" to represent DSA preconditioned GMRES using zero initial guess and ``PGMRES-ROMIG" to represent DSA preconditioned GMRES using ROMIG.}
Recall that we denote the tolerance in the stopping criteria of SISA as $\epssisa$, the truncation tolerance to  POD basis in \eqref{eq:pod_truncation} as $\epspod$, and the window size to assemble the snapshot matrix for the correction equation as $\mfw$.
We define the average number of \pzc{transport sweeps} for the test set as 
\begin{equation*}
    \pzc{n_{\textrm{sweep}}} = \frac{\sum_{\bmu\in\mathcal{P}_{\textrm{train}}}\textrm{number of \pzc{transport sweeps} for the convergence for the parameter}\bmu }{\textrm{total number of parameters in }\mathcal{P}_{\textrm{test}}}.
\end{equation*}
\pzc{In high dimensions, the average computational time is approximately proportional to $n_{\textrm{sweep}}$ due to many angular directions sampled in the angular discretization.}
Similarly, we define the average relative computational time with respect to SI-DSA as $T_{\textrm{rel}}$.  
To verify that \pzc{iterative solvers converge to the correct solution}, we \pzc{compute the average residual of equation \eqref{eq:rho_matrix_vec_def}:} 
 \begin{equation*}
     \pzc{\mathcal{R}_{\infty}^\rho = \frac{\sum_{\bmu\in\mathcal{P}_{\textrm{train}}}||(\mathbf{I}-\mathbf{L})\brho_{\bmu}-\overline{\bb}_{\bmu}||_\infty}{\textrm{total number of parameters in }\mathcal{P}_{\textrm{test}}}},
 \end{equation*}
where \pzc{$\brho_{\bmu}$ is the numerical solution obtained by the underlying iterative solver}.

\begin{table}[htbp]
\centering
  \begin{tabular}{|l|c|c|c|c|c|c|c|c|c|c|c|}
    \hline
\pzc{DSA} & \pzc{Source Iteration with Diffusion Synthetic Acceleration and $0$ initial guess}\\ \hline
\pzc{SI-ROMIG} & \pzc{SI-DSA  with ROM based initial guess}\\ \hline
\pzc{ROMSA-$\mfw$} & \pzc{ROMSA using a ROM generated with window size $\mfw$}\\ \hline
\multirow{2}{7em}{\pzc{ROMSAD-$\mfw,\theta$}}
  & \pzc{ROMSAD using a ROM generated with window size $\mfw$,} \\ 
  &\pzc{and ROMSA is allowed to be applied in at most the first $\theta$ iterations} \\ \hline
\pzc{PGMRES} & \pzc{DSA preconditioned GMRES using $0$ initial guess}\\ \hline
\pzc{PGMRES-ROMIG} & \pzc{DSA preconditioned GMRES using ROM based initial guess}\\ \hline
{\pzc{$n_{\textrm{sweep}}$}} & \pzc{Average number of transport sweep required by convergence}\\ &\pzc{for parameters in the test set} \\ \hline
{\pzc{$T_{\textrm{rel}}$}} & \pzc{Relative computational time with respect to SI-DSA using $0$ initial guess}\\ \hline
\pzc{$\epspod$} & \pzc{Truncation tolerance to generate POD basis in \eqref{eq:pod_truncation}} \\ \hline
\pzc{$\epssisa$} & \pzc{Tolerance in the stopping criteria of SISA}\\ \hline
\pzc{$\mathcal{R}_{\infty}^\rho$} & \pzc{Average residual of equation \eqref{eq:rho_matrix_vec_def}}\\ \hline
 \end{tabular}
 \caption{Notations used in Sec. \ref{sec:numerical}. \label{tab:notation}}
\end{table}

\subsection{1D slab geometry}
In our tests for 1D slab geometry,  $16$ Gauss-Legendre quadrature points for $[-1,1]$ is used for angular discretization. \pzc{The main purpose of the 1D tests is to investigate the influence of the POD truncation tolerance $\epspod$ and the window size $\mfw$ on the performance of proposed methods.}

\pzc{In 1D tests, we omit the relative computational time for the following reasons.  Since we only use $16$ quadrature points in our angular discretization, the computational time for solving the diffusion correction equation in DSA is not negligible compared to the transport sweep. On the other hand, the computational time for solving the small size reduced order correction equation in ROMSA is negligible compared to the transport sweep. However, in more practical high dimensional problems, a large number of angular directions is required, resulting in the correction step's computational time being much smaller than the transport sweep, regardless of whether DSA, ROMSA, or ROMSAD is used. Consequently, computational savings gained by ROMSA or ROMSAD are more pronounced in 1D.}
\subsubsection{Cross-regime problem\label{sec:num-1d-cross}}
We consider a parametric problem with a spatially varying scattering cross section. The set-up of the problem is as follows:
\begin{align*}
\Omega_{\bx}=[0,10] ,\;G = 0.01, \;
\sigma_a = 0,\;
\sigma_s(x) = 0.1+\mu_s x,\;
 \;\sigma_t = \sigma_a+\sigma_s, \;
f(0,\xi) = 0\;\text{with } \xi>0, \;f(10,\xi) = 0 \;\text{with } \xi< 0.  
\end{align*}
The slope of the linear scattering cross section $\mu_s\in[10,20]$ is the parameter of this problem. 
We use a uniform mesh with $N_{\bx}=400$ to partition the computational domain. The tolerance in the stopping criteria of SISA is set to  $\epssisa=10^{-11}$.

As we move from the left to the right of the computational domain, the scattering cross section $\sigma_s(x)$ continuously changes from $0.1$ to $0.1+10\mu_s$ which is bigger than $100$ for $\mu_s \in[10,20]$. In other words, there is a smooth transition in the material property from transport dominance to scattering dominance. The slope $\mu_s$ determines how fast this transition is. 

\begin{figure}[]
  \begin{center} 
\includegraphics[width=0.45\textwidth]{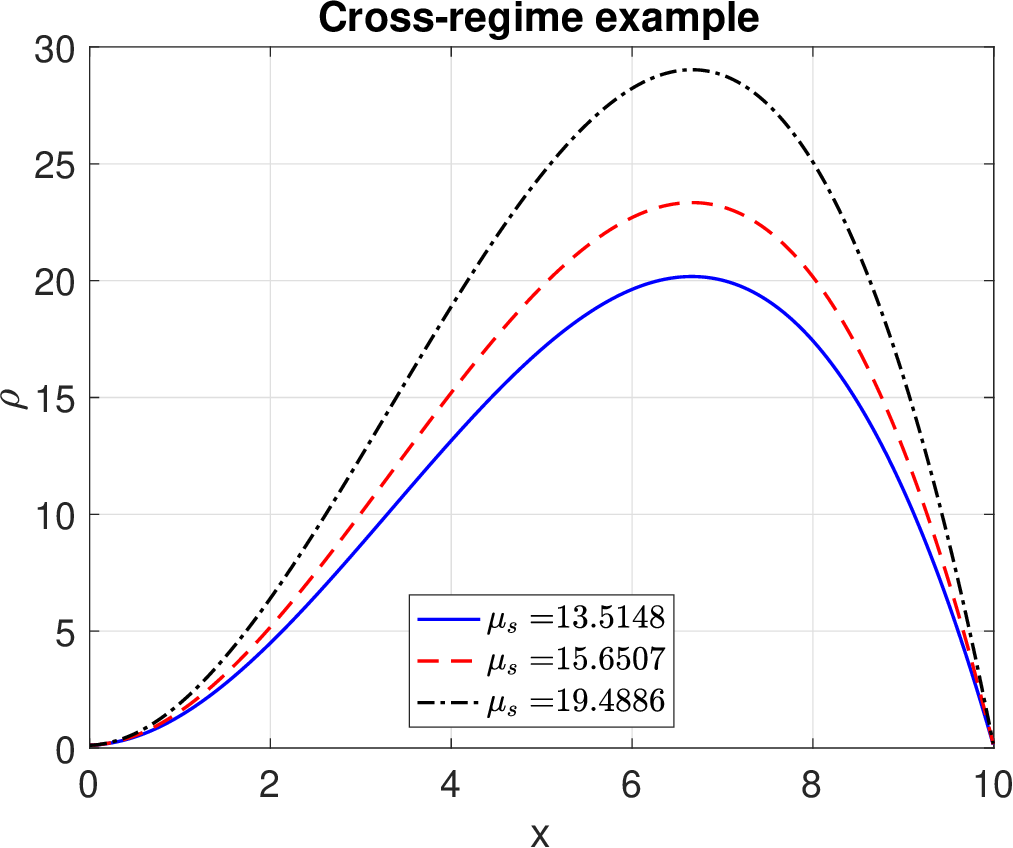}
\includegraphics[width=0.45\textwidth]{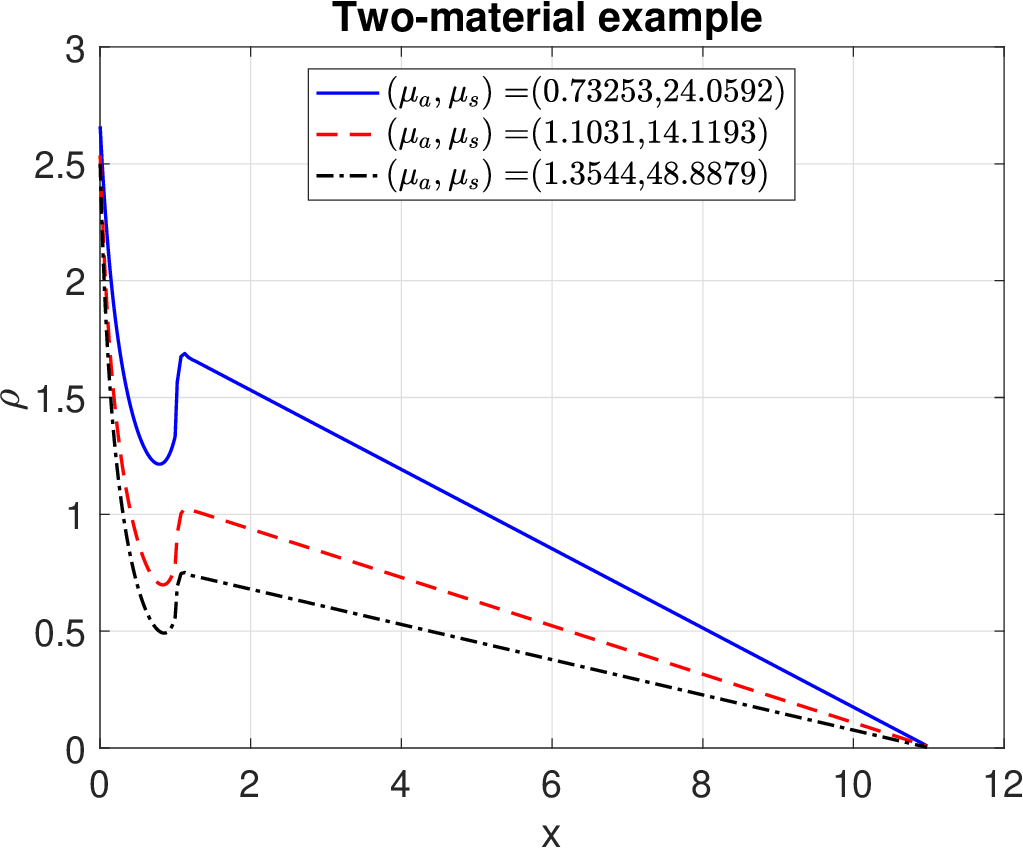}
  \caption{Reference solutions for $3$ samples in the test set. Left: the cross-regime example in Sec. \ref{sec:num-1d-cross}. Right: the two-material example in Sec. \ref{sec:num-1d-two-material}. \label{fig:1d-reference}}
  \end{center}
\end{figure}

The training set $\mathcal{P}_{\textrm{train}}=\{\mu_{s,j}=10+j\Delta\mu,\;\Delta\mu=\frac{1}{4},\;j=0,\dots,40\}$ consists of $41$ uniformly distributed samples in $[10,20]$. To test the performance of the proposed methods, we randomly sample $20$ different values of $\mu_s$ from $[10,20]$. Solutions 
corresponding to $3$ test samples are presented in the left picture of Fig. \ref{fig:1d-reference}.
Using this cross-regime problem, we investigate the influence of the POD truncation tolerance $\epspod$ and the window size $\mfw$ on the performance of \pzc{ROMIG}, ROMSA and ROMSAD. 

Dimensions of the reduced order spaces in our ROMs for the parametric problem ($r_p$) and the corresponding kinetic correction equations ($r_c$) are displayed 
in Tab. \ref{tab:1d_cross_regime_rank}. With the same truncation tolerance $\epspod$, We observe that $r_p$ is always not greater than $r_c$. In addition, for the same tolerance $\epspod$, the dimension of the reduced order space for the correction equation, $r_c$, grows as the window size $\mfw$ increases.

In Tab. \ref{tab:1d_cross_regime_performance}, we compare the performance of various methods.  Key observations are as follows.
\begin{enumerate}
    \item \textbf{Influence of the POD truncation tolerance $\epspod$:} For this problem, the performance of \pzc{ROMIG} is insensitive to the choice of the tolerance $\epspod$, and it converges with nearly $60\%$ of iterations with respect to DSA.

    ROMSA is sensitive to the choice of the tolerance $\epspod$ for this problem. With $\epspod=10^{-9}$, ROMSA needs even more iterations for convergence than DSA,  regardless of the window size $\mfw$.  The results for ROMSA-$1$ are identical with $\epspod=10^{-9}$ and $\epspod=10^{-10}$, as the underlying reduced order spaces are the same. 
    When $\epspod=10^{-11}$, ROMSA-$1$ converges with $47\%$ number of iterations with respect to DSA. 
    For $\epspod=10^{-10}$ and $\epspod=10^{-11}$, ROMSA-$3$ \pzc{converges with fewer transport sweeps} than SI-DSA, regardless of whether the zero or the ROM-based initial guess is used.

    ROMSAD is more robust than ROMSA concerning the choice of the tolerance $\epspod$. Compared with DSA and \pzc{ROMIG}, ROMSAD consistently \pzc{converges with fewer transport sweeps}, regardless of the window size $\mfw$. Due to its robustness, ROMSAD-$\mfw,3$ is significantly more efficient than ROMSA-$\mfw$ for $\epspod=10^{-9}$ and $\epspod=10^{-10}$. When $\epspod=10^{-11}$, the underlying ROM for the correction equation becomes accurate, and ROMSAD-$\mfw,3$ is slightly more efficient than ROMSA-$\mfw$.
    
    \item \textbf{Influence of $\mfw$ on ROMSAD:} When $\epspod=10^{-9}$, ROMSAD-$\mfw,3$ converges with a comparable number of iterations for window sizes \pzc{$\mfw=1,3$}.
    When $\epspod=10^{-10}$ or $10^{-11}$, ROMSAD converges faster if the underlying ROM is generated with a wider window size $\mfw$.  

\end{enumerate}

For this problem, ROMSAD and \pzc{ROMIG} are more robust than ROMSA, and more efficient than DSA. In the rest of the numerical section, we will mainly focus on comparing these two methods with DSA.
\begin{table}[htbp]
\centering
  \begin{tabular}{|l|c|c|c|c|c|c|c|c|c|c|c|}
    \hline
     $\epspod$ & $r_p$ & $r_c$ with $\mfw=1$ &  $r_c$ with $\mfw=3$  
 \\ \hline	
$10^{-9}$ & $6$ & $6$ & $15$ 
\\ \hline
$10^{-10}$ & $6$ &$6$ & $18$ 
\\ \hline
$10^{-11}$ & $7$ & $7$ & $21$ 
\\ \hline
 \end{tabular}
 \caption{ Dimensions of the reduced order spaces for the parametric problem ($r_p$) and the correction equation ($r_c$) for the cross-regime problem in Sec. \ref{sec:num-1d-cross} with different POD truncation tolerances $\epspod$. \label{tab:1d_cross_regime_rank}}
 \end{table}

\begin{table}[htbp]
\centering
\begin{subtable}{\textwidth}
\centering
   \begin{tabular}{|l|c|c|c|c|c|c|c|c|c|c|c|}
    \hline
    & DSA & \pzc{ROMIG} & \pzc{ROMSA}-$1$ & \pzc{ROMSAD}-$1,3$ & \pzc{ROMSA}-$3$ & \pzc{ROMSAD}-$3,3$ 
    \\ \hline   
    $\pzc{{n}_{\textrm{sweep}}}$ & 15 &  9 & 23.75 & 6.55 & 42.25 & 6.85 
    \\ \hline
    $\pzc{\mathcal{R}^\rho_{\infty}}$ & \pzc{1.02e-12}  & \pzc{2.33e-12} & \pzc{8.85e-12} & \pzc{1.85e-12} & \pzc{8.98e-12} & \pzc{2.82e-12} 
    \\ \hline
 \end{tabular}
     \caption{$\epspod=10^{-9}$, $\epssisa=10^{-11}$\label{tab:1d_cross_regime_a}}
 \end{subtable}\\
 \medskip
\begin{subtable}{\textwidth}
\centering
  \begin{tabular}{|l|c|c|c|c|c|c|c|c|c|c|c|}
    \hline
    & DSA & \pzc{ROMIG} & \pzc{ROMSA-$1$} & \pzc{ROMSAD}-$1,3$ & \pzc{ROMSA}-$3$ & \pzc{ROMSAD}-$3,3$ 
    \\ \hline 
    $\pzc{{n}_{\textrm{sweep}}}$ & 15 &  9 & 23.75 & 6.6 & 12.05 & 5.8 
    \\ \hline
    $\pzc{\mathcal{R}^\rho_{\infty}}$ & \pzc{1.02e-12}  & \pzc{2.33e-12} & \pzc{8.85e-12} & \pzc{1.77e-12} & \pzc{6.04e-12} & \pzc{1.41e-12} 
    \\ \hline
 \end{tabular}
 \caption{$\epspod=10^{-10}$, $\epssisa=10^{-11}$\label{tab:1d_cross_regime_b}}
 \end{subtable}\\
 \medskip
\begin{subtable}{\textwidth}
\centering
 \begin{tabular}{|l|c|c|c|c|c|c|c|c|c|c|c|}
    \hline
    & DSA & \pzc{ROMIG} & \pzc{ROMSA}-$1$ & \pzc{ROMSAD}-$1,3$ & \pzc{ROMSA}-$3$ & \pzc{ROMSAD}-$3,3$ 
    \\ \hline  
    $\pzc{{n}_{\textrm{sweep}}}$ & 15 &  9 & 7.05 & 5.65 & 4.55 & 4.5 
    \\ \hline
    ${\mathcal{R}}^{\rho}_{\infty}$ & \pzc{1.02e-12} &  \pzc{2.93e-12} & \pzc{5.82e-12} & \pzc{8.67e-13} & \pzc{3.27e-12} & \pzc{1.59e-12} 
    \\ \hline
 \end{tabular}
 \caption{$\epspod=10^{-11}$, $\epssisa=10^{-11}$\label{tab:1d_cross_regime_c}}
 \end{subtable}
 \caption{Results for the cross-regime problem in 1D slab geometry (see Sec. \ref{sec:num-1d-cross}). \label{tab:1d_cross_regime_performance}}
 \end{table}

\subsubsection{Two-material problem\label{sec:num-1d-two-material}} 
We consider a parametric two-material problem with the following set-up:
\begin{align*}
&\Omega_{\bx}=[0,11] ,\;G = 0, \;
\sigma_a(x) = \begin{cases}
            \mu_a,\;0<x\leq 1,\\
		0,\; 1<x<11,
		\end{cases}\;
\sigma_s(x) = \begin{cases}
	0,\;0<x\leq1,\\
		\mu_s,\;  1<x<11,
		\end{cases}
 \;\sigma_t(x) = \sigma_a(x)+\sigma_s(x), \notag\\
&f(0,\xi) = 5\;\text{with } \xi>0, \;f(11,\xi) = 0 \;\text{with } \xi< 0.  
\end{align*}
There is a pure absorption region ($[0,1]$) and a pure scattering region ($[1,11]$) in the computational domain.
Parameters $\mu_a\in[0.5,1.5]$ and $\mu_s\in[10,50]$ are the strength of the absorption and scattering cross sections in the absorption and scattering regions, respectively.
A non-uniform mesh with $\Delta x=\Delta x_1=\frac{1}{100}$ on $[0,1]$ and $\Delta x=\Delta x_2 =\frac{1}{10}$ on $[1,11]$ is used to partition the computational domain. The tolerance in the stopping criteria of SISA is set to $\epssisa=10^{-12}$.

The training set for this problem is
\[\mathcal{P}_{\textrm{train}}=\{(\mu_{a,m},\mu_{s,n}):\;\mu_{a,m}=0.5+0.1m,\;\mu_{s,n}=10+n,0\leq m\leq 10,0\leq n\leq 40\}.\]
To test the performance of the proposed methods, we randomly sample $20$ pairs of $(\mu_a,\mu_s)$ from $[0.5,1.5]\times[10,50]$. Solutions 
corresponding to $3$ samples in the test set are presented in the right picture of Fig. \ref{fig:1d-reference}.

The results are summarized in Tab. \ref{tab:1d_two_material}. 
With $\theta=1,2,3$, i.e. the maximum number of iterations allowed to use ROMSA in ROMSAD, ROMSAD-$\mfw,\theta$ turn out to have identical performance for this problem, when the same window size $\mfw$ ($\mfw=3,5$) is applied. With the POD truncation tolerances $\epspod=10^{-6}$, $10^{-9}$ and $10^{-10}$, both \pzc{ROMIG} and ROMSAD leads to convergence faster than DSA. The efficiency of \pzc{ROMIG} and ROMSAD improves, as we decrease the tolerance $\epspod$.
For this example, \pzc{ROMIG} is slightly more efficient than ROMSAD, when their underlying ROMs are generated with the same tolerance $\epspod$. 

\begin{table}[htbp]
\centering
   \begin{tabular}{|l|c|c|c|c|c|c|c|c|c|c|c|}
    \hline
&  $\epspod$ & DSA & \pzc{ROMIG} & \pzc{ROMSAD}-$3,\theta$ & \pzc{ROMSAD}-$5,\theta$ \\ \hline
 $\pzc{n_{\textrm{sweep}}}$ & \multirow{3}{*}{$10^{-6}$} &  14.6 & 8.1 & 8.9 & 9.05   \\ \cline{1-1} \cline{3-6}
 $\pzc{\mathcal{R}^\rho_{\infty}}$ &  & \pzc{1.94e-13} &  \pzc{1.46e-13} & \pzc{1.80e-13} & \pzc{2.05e-13} \\ \hline
 $\pzc{n_{\textrm{sweep}}}$ & \multirow{3}{*}{$10^{-9}$} &  14.6 & 4.2 & 5.35 & 5.35   \\ \cline{1-1} \cline{3-6}
 $\pzc{\mathcal{R}^\rho_{\infty}}$ &  & \pzc{1.94e-13} & \pzc{1.48e-13} & \pzc{2.15e-13} & \pzc{2.24e-13} \\ \hline
  $\pzc{n_{\textrm{sweep}}}$ & \multirow{3}{*}{$10^{-10}$} &  14.6 & 3.35 & 4.25 & 4.25   \\ \cline{1-1} \cline{3-6}
 $\pzc{\mathcal{R}^\rho_{\infty}}$ &  & \pzc{1.94e-13}  & \pzc{1.01e-13} & \pzc{2.35e-13} & \pzc{1.91e-13} \\ \hline
 \end{tabular}
     \caption{Results for the two-material problem in Sec. \ref{sec:num-1d-two-material}. The most number of iterations applying ROMSA in ROMSAD is $\theta=1,3,5$. Performance of ROMSAD-$3,\theta$ and ROMSAD-$5,\theta$ are identical to each other. Dimensions of the reduced order spaces: when $\epspod=10^{-6}$, $r_p=21$, $r_c=49$ for $\mfw=3$ and $r_c=54$ for $\mfw=5$; when $\epspod=10^{-9}$, $r_p=44$, $r_c=101$ for $\mfw=3$ and $r_c=120$ for $\mfw=5$;  when $\epspod=10^{-10}$, $r_p=53$, $r_c=123$ for $\mfw=3$ and $r_c=147$ for $\mfw=5$.\label{tab:1d_two_material}}
\end{table}

\subsection{2D $X$-$Y$ geometry and $\mathbb{S}^2$ as angular space}
To \pzc{compare} the performance of \pzc{different} methods, we consider a \pzc{seris} of problems in 2D $X$-$Y$ geometry with $\mathbb{S}^2$ as the angular space. For the $S_N$ angular discretization, we use the CL($N_\phi$,$N_{v_z}$) quadrature rule, which is the tensor product of the normalized $N_\phi$-points Chebyshev quadrature rule for the unit circle and the normalized $N_{v_z}$-points Gauss-Legendre quadrature rule for $v_z\in[-1,1]$.  

\pzc{In the 2D X-Y geometry, many angular directions are sampled. Hence, the computational time of transport sweeps are more dominant compared to the Synthetic Acceleration step. As a result, the computational time is nearly proportional to the number of transport sweeps required for convergence.}
\subsubsection{Homogeneous medium \label{sec:line-source}}
We consider a parametric problem with zero inflow boundary conditions and a homogeneous background medium:
\begin{equation*}
\Omega_{\bx}=[0,1]^2 ,\;G(x,y) = \exp(-100( (x-0.5)^2+(y-0.5)^2)), \;
\sigma_a(\bx) = 0,\;
\sigma_s(\bx) = \mu_s,\;\mu_s\in[0.9,1.1].
\end{equation*}
We use an $80\times 80$ uniform rectangular mesh to partition $\Omega_{\bx}$ and the CL-($30,6$) quadrature rule  to discretize the angular space. The tolerance in the stopping criteria of SISA is set to $\epssisa=10^{-12}$. 
\pzc{The tolerance of relative residual in PGMRES is $10^{-11}$. As shown in Tab. \ref{tab:line-source}, converged solutions given by SI and PGMRES result in residuals with comparable magnitude.}

\begin{figure}[]
  \begin{center} 
\includegraphics[width=0.45\textwidth]{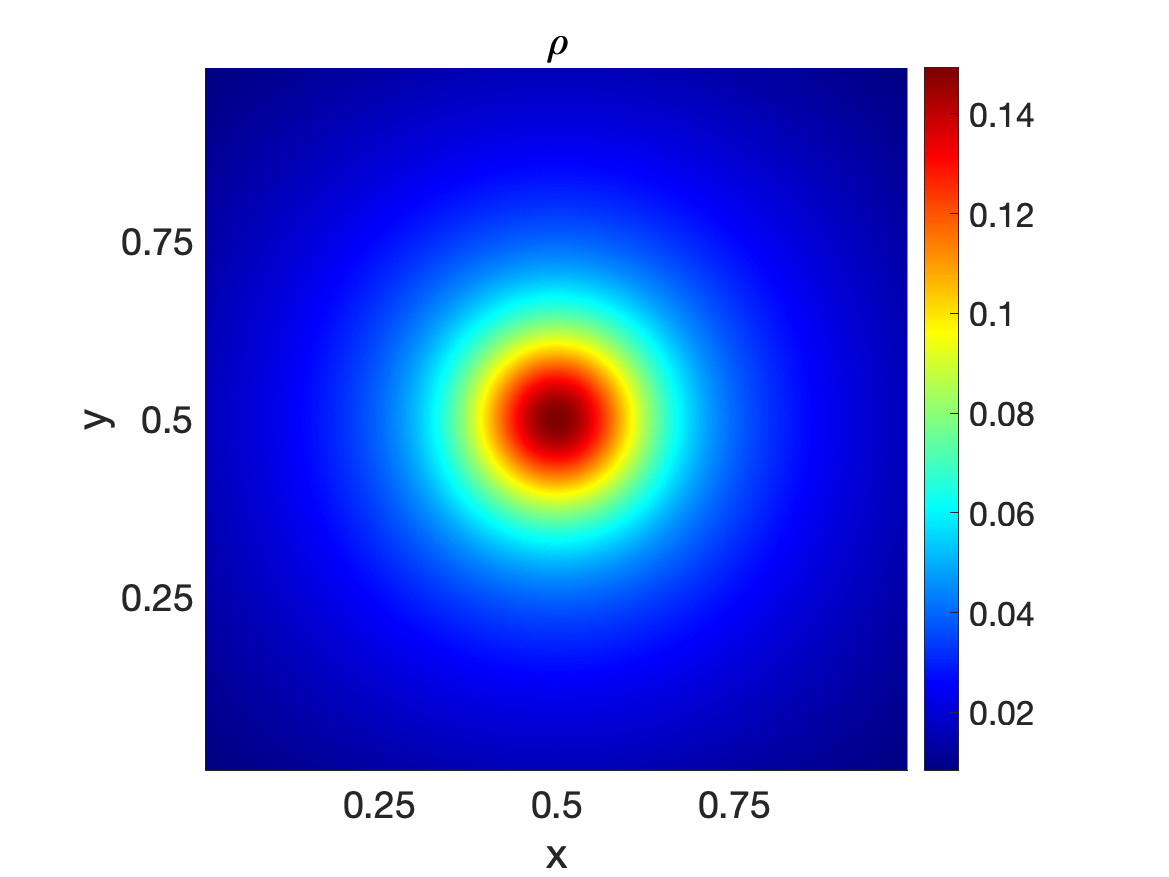}
\includegraphics[width=0.45\textwidth]{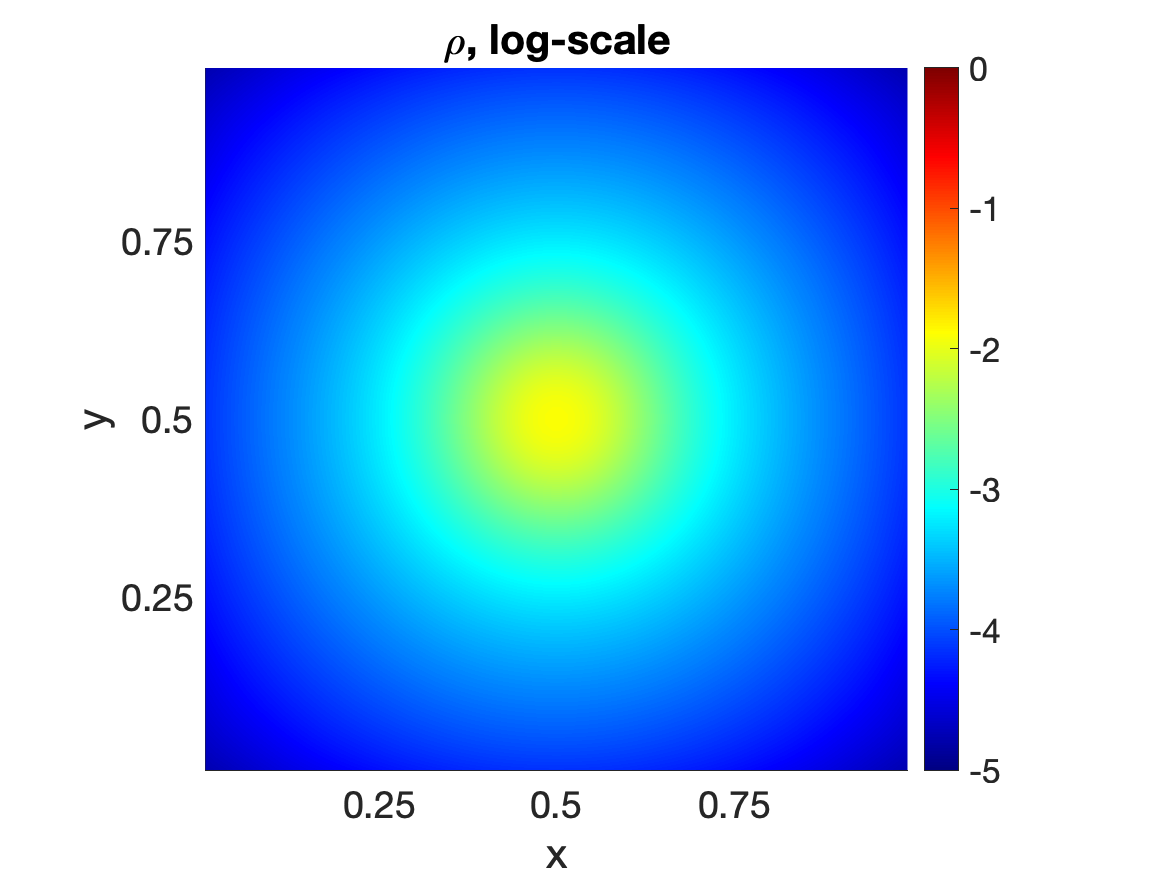}
  \caption{The reference solution for the homogeneous medium problem in Sec. \ref{sec:line-source} with $\mu_s=1.0101$. The right figure is the  solution under log-scale. \label{fig:line-source}}
  \end{center}
\end{figure}

For this problem, we construct our ROMs using a uniformly sampled training set with $21$ samples: $\mathcal{P}_{\textrm{train}}=\{\mu_{s,m}=0.9+m\Delta\mu,\;\Delta\mu=\frac{1}{200},0\leq m\leq 20\}$. We test the performance of the proposed methods with $10$ randomly sampled values of $\mu_s$ from $[0.9,1.1]$. POD basis is generated with the tolerance $\epspod=10^{-9}$. The solution for the test sample $\mu_s=1.0101$ is presented in Figure \ref{fig:line-source}. 

\textbf{\pzc{Online efficiency:}} Results for the test set are summarized in Tab. \ref{tab:line-source}.  \pzc{Under the SI framework, ROMIG, ROMSA and ROMSAD yield
slightly more than $5$ times acceleration compared to DSA}.  \pzc{SI-ROMSAD is approximately $4$ times as efficient as DSA preconditioned GMRES} \pzc{ROMIG is able to accelerate DSA preconditioned GMRES. However, it is still slightly slower than SI-DSA using ROMIG or SI-ROMSAD due to the extra transport sweeps discussed in  Sec. \ref{sec:details}. }

\textbf{\pzc{Offline costs:}} \pzc{The offline cost to construct ROMs includes computing numerical solutions for $21$ training samples and constructing reduced order basis and operators.  In Tab. \ref{tab:offline-line-source}, we summarize the relative offline computational time needed to construct reduced order basis and operators  with respect to the average computational time of solving RTE for parameters in the test set with SI-DSA. For the parametric problem, the relative computational time for constructing the reduced order bases and operators is approximately $0.711\%$ and $0.238\%$, respectively. For the correction equation, the relative time for constructing the reduced order basis and operators is approximately $2.722\%$ and $0.471\%$, respectively.} 
\begin{table}
\centering
 \begin{tabular}{|l|c|c|c|c|c|c|c|c|c|c|c|}
    \hline
    & DSA  & \pzc{ROMIG} & \pzc{ROMSA}-$3$ & \pzc{ROMSAD}-$3,3$  & \pzc{PGMRES} & \pzc{PGMRES-ROMIG} \\ \hline  
    $\pzc{{n}_{\textrm{sweep}}}$ & 14.6 &  3.0 & 2.9 & 2.9 & \pzc{11} & \pzc{4.8}  \\ \hline
    ${T}_{\textrm{rel}}$ & 100\%  & \pzc{19.512\%} & \pzc{16.578\%} & \pzc{17.853\%} & \pzc{78.778\%} & \pzc{33.966\%} \\ \hline
    $\pzc{\mathcal{R}^\rho_{\infty}}$ & \pzc{1.05e-13} &  \pzc{1.44e-13} & \pzc{1.26e-13} & \pzc{1.22e-13} & \pzc{5.45e-14} & \pzc{1.14e-13}   \\ \hline
 \end{tabular}
 \caption{Results for the problem with homogeneous medium in Sec. \ref{sec:line-source}. Dimensions of reduced order spaces: $r_p=5$ and $r_c=14$.\label{tab:line-source}}
 \end{table}
\begin{table}
\centering
 \begin{tabular}{|l|c|c|c|c|c|c|c|c|c|c|c|}
    \hline
    & \pzc{Reduced basis construction} & \pzc{Reduced operator construction}   \\ \hline  
\pzc{Parametric problem} &   \pzc{$0.711\%$} &    \pzc{$0.238\%$}\\ \hline
\pzc{Correction equation} &  \pzc{$2.722\%$} & \pzc{$0.471\%$}  \\ \hline
 \end{tabular}
 \caption{\pzc{Relative offline computational cost of constructing reduced basis and reduced order operators for the test in Sec. \ref{sec:line-source}. Dimensions of reduced order spaces: $r_p=5$ and $r_c=14$.}\label{tab:offline-line-source}}
 \end{table}
\subsubsection{Lattice problem \label{sec:lattice}}
\begin{figure}[]
  \begin{center} 
\includegraphics[width=0.45\textwidth]{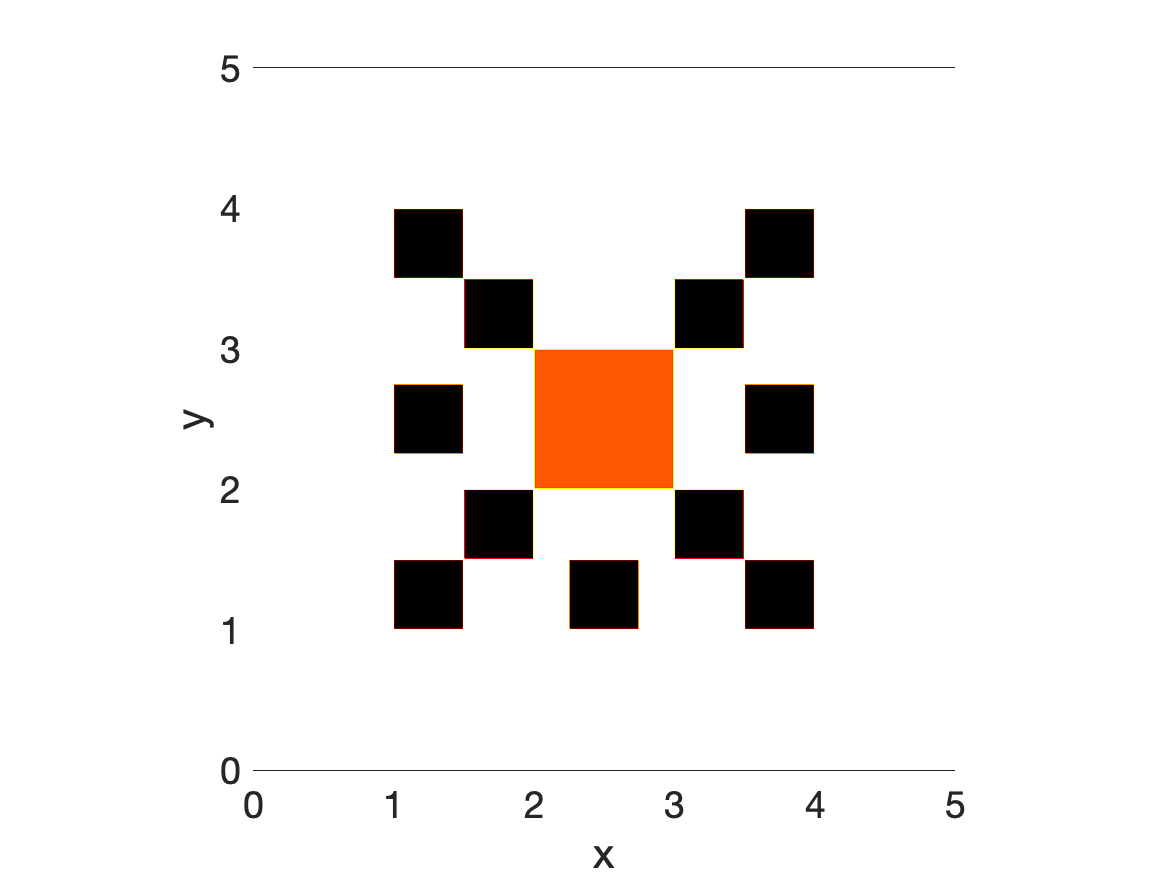}
\includegraphics[width=0.45\textwidth]{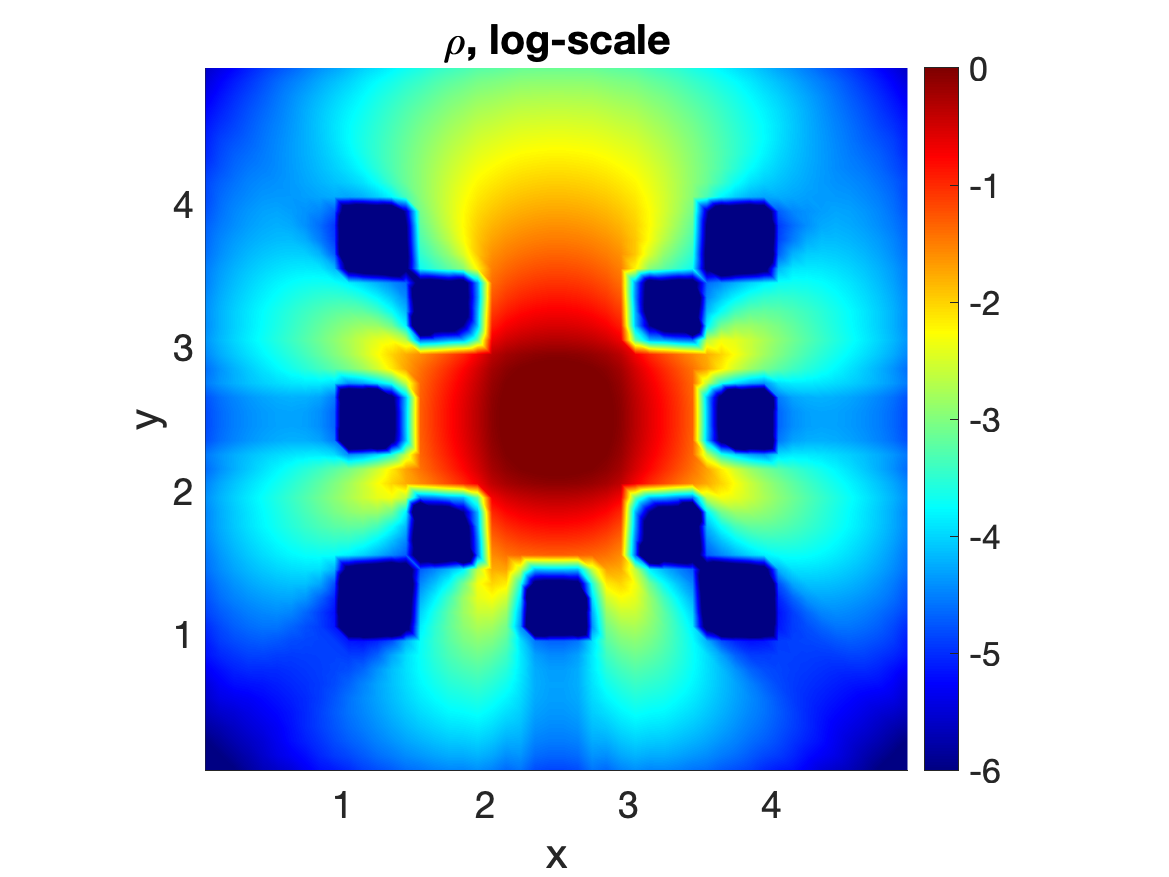}
  \caption{The set-up and a reference solution for the lattice problem in Sec. \ref{sec:lattice}. Left: the set-up for the lattice problem. Black: pure absorption regions with $(\sigma_a,\sigma_s)=(\mu_a,0)$. White and orange: pure scattering regions with $(\sigma_a,\sigma_s)=(0,\mu_s)$. Orange: constant source term with $G=1$. Right: the reference solution under log-scale for $(\mu_a,\mu_s)=(97.235,1.051)$. \label{fig:lattice}}
  \end{center}
\end{figure}

We consider a parametric lattice problem with zero inflow boundary conditions in the computational domain $\Omega_{\bx}=[0,5]^2$. The geometry set-up is presented in the left picture of Fig. \ref{fig:lattice}. In Fig. \ref{fig:lattice}, black regions are pure absorption regions  with $(\sigma_a,\sigma_s)=(\mu_a,0)$. Other regions are pure scattering regions with $(\sigma_a,\sigma_s)=(0,\mu_s)$. The parameter $\bmu=(\mu_a,\mu_s)\in[95,105]\times[0.5,1.5]$ determines the strength of absorption in the pure absorption region and the strength of scattering in the pure scattering region. The source term 
\begin{align*}
G(x,y) = \begin{cases}
         1.0,\quad\text{if}\quad |x-2.5|<0.5\;\text{and}\quad|y-2.5|<0.5,\\
         0, \quad\text{otherwise},
         \end{cases}
\end{align*}
is imposed in the orange region.
We use a uniform mesh with $50\times 50$ elements to partition the computational domain and discretize the angular space with the CL($40,6$) quadrature rule. The tolerance for the SISA stopping criteria is set to $\epssisa=10^{-12}$. \pzc{The tolerance of relative residual in PGMRES is $10^{-12}$. As shown in Tab. \ref{tab:lattice}, converged solutions given by SI and PGMRES result in residuals of comparable magnitude.}

The training set $$\mathcal{P}_{\textrm{train}}=\{(\mu_a,\mu_s)=(95+i\Delta\mu_a,0.5+0.1j\Delta\mu_s),\Delta\mu_a=1,\Delta \mu_s=0.1,i,j=0,1,\dots,10\}$$
has $121$ pairs of uniformly sampled $(\mu_a,\mu_s)$ in total. We test the performance of our methods with $10$ randomly sampled $(\mu_a,\mu_s)$ from $[95,105]\times[0.5,1.5]$.  A reference solution corresponding to a test sample is presented in the right picture of Fig. \ref{fig:lattice}.

\textbf{\pzc{Online efficiency:}} Results for parameters in the test set are summarized in Tab. \ref{tab:lattice}. Various POD truncation tolerances are applied to generate ROMs.  \pzc{When $\epspod=10^{-9}$}, \pzc{SI-DSA using ROMIG} is approximately $1.6$ times as fast as ROMSAD-$3,5$, \pzc{$2.6$ times as fast as DSA preconditioned GMRES} and $4$ times as fast as SI-DSA \pzc{using the zero initial guess}. 
When $\epspod=10^{-10}$, \pzc{SI-DSA using ROMIG} and ROMSAD-$3,5$ have comparable speeds, leading to a nearly $5$ times acceleration compared to \pzc{SI-DSA and slightly more than $3.4$ times acceleration compared to DSA preconditioned GMRES}. 
When $\epspod=10^{-11}$, ROMSAD-$3,5$ becomes the fastest, leading to a convergence with approximately $11.41\%$ computational time and $12.69\%$ number of iterations compared to \pzc{SI-DSA}. 
\pzc{Still, due to the extra transport sweeps discussed in Sec. \ref{sec:details}, DSA preconditioned GMRES using ROMIG is slightly less efficient than SI-DSA using ROMIG, but it is at least  $2$ times as fast as DSA preconditioned GMRES using the zero initial guess.}
Overall, for this example, as the POD truncation tolerance $\epspod$ decreases, both of \pzc{ROMIG} and ROMSAD-$3,5$ gain performance boosts, while the boost gained by ROMSAD-$3,5$ is larger. 

\pzc{\textbf{Offline costs:} The main offline cost is to compute numerical solutions for $121$ training samples.  In Tab. \ref{tab:offline-lattice}, we present the relative offline computational time needed to construct reduced order basis and operators  with respect to the average computational time of solving RTE for parameters in the test set with SI-DSA. The computational time of constructing reduced basis are almost identical for various $\epspod$, since SVD of the full snapshot matrix takes up the most computational time. The relative time to construct reduced basis is approximately $10\%$ for the parametric problem itself and $84\%$ for the correction equation. The relative time to construct reduced order operators increases as the dimension of the reduced order space grows and is less than $1.7\%$ for the parametric problem itself and $6.0\%$ for the correction equation.} 
\begin{table}[htbp]
\centering
   \begin{tabular}{|l|c|c|c|c|c|c|c|c|c|c|c|}
    \hline
&  $\epspod$ & DSA & \pzc{ROMIG} & \pzc{ROMSAD}-$3,5$ & \pzc{PGMRES} & \pzc{PGMRES-ROMIG}  \\ \hline
 $\pzc{n_{\textrm{sweep}}}$ & \multirow{3}{*}{$10^{-9}$} &  19.7& 5.1 & 8.0 & \pzc{13} & \pzc{5.9}  \\ \cline{1-1} \cline{3-7}
 $T_{\textrm{rel}}$ &  & 100\% &  \pzc{25.493\%} & \pzc{39.739\%} & \pzc{62.442\%} & \pzc{28.083\%}  \\ \cline{1-1} \cline{3-7}
 $\pzc{\mathcal{R}^\rho_{\infty}}$ &  & \pzc{2.33e-13} & \pzc{1.94e-13} & \pzc{4.78e-13} & \pzc{2.57e-13} & \pzc{1.26e-13}  \\ \hline
 $\pzc{n_{\textrm{sweep}}}$ & \multirow{3}{*}{$10^{-10}$} &  19.7 & 3.9 & 4.2 & \pzc{13} & \pzc{5.1}  \\ \cline{1-1} \cline{3-7}
 $T_{\textrm{rel}}$ &  & 100\%  & \pzc{19.336\%} & \pzc{20.318\%} & \pzc{62.442\%}  & \pzc{24.735\%}  \\ \cline{1-1} \cline{3-7}
 $\pzc{\mathcal{R}^\rho_{\infty}}$ &  & \pzc{2.33e-13}  & \pzc{1.66e-13} & \pzc{2.58e-13} & \pzc{2.57e-13} & \pzc{2.53e-13} \\ \hline
 $\pzc{n_{\textrm{sweep}}}$ & \multirow{3}{*}{$10^{-11}$} &  19.7 &  3.5 & 2.5 & \pzc{13} & \pzc{4.4}  \\ \cline{1-1} \cline{3-7}
  $T_{\textrm{rel}}$ &  & 100\% &  \pzc{17.316\%} & \pzc{11.387\%}    & \pzc{62.442\%}  & \pzc{21.408}\% \\ \cline{1-1} \cline{3-7}
 $\pzc{\mathcal{R}^\rho_{\infty}}$ &  & \pzc{2.33e-13} & \pzc{1.62e-13} & \pzc{9.72e-14} & \pzc{2.57e-13} & \pzc{5.02e-13}\\ \hline
 \end{tabular}
     \caption{Results for the lattice problem in Sec. \ref{sec:lattice}. Dimensions of the reduced order spaces: $r_p=25$ and $r_c=66$ for $\epspod=10^{-9}$; $r_p=29$ and $r_c=80$ for $\epspod=10^{-10}$;  $r_p=34$ and $r_c=96$ for $\epspod=10^{-11}$.\label{tab:lattice}}
\end{table}
\begin{table}
\centering
 \begin{tabular}{|l|c|p{22mm}|p{28mm}|p{23mm}|}
    \hline
    & \pzc{$\epspod$} & \pzc{Reduced basis construction} & \pzc{Reduced operator construction} & \pzc{Reduced space dimension}  \\ \hline  
\multirow{3}{2cm}{\pzc{Parametric problem}} & \pzc{$10^{-9}$} &   \pzc{$10.360\%$} &    \pzc{$1.355\%$} & \pzc{$25$}\\ \cline{2-5}
& \pzc{$10^{-10}$} &   \pzc{$10.806\%$} &    \pzc{$1.372\%$} & \pzc{$29$} \\ \cline{2-5}
& \pzc{$10^{-11}$} &   \pzc{$9.322\%$} &    \pzc{$1.643\%$} &  \pzc{$34$}\\ \hline
\multirow{3}{2cm}{\pzc{Correction equation}}  &  \pzc{$10^{-9}$} &\pzc{$83.584\%$} & \pzc{$3.394\%$} &  \pzc{$66$} \\ \cline{2-5}
&\pzc{$10^{-10}$} &\pzc{$84.784\%$} & \pzc{$4.067\%$}  &  \pzc{$80$}\\ \cline{2-5}
&\pzc{$10^{-11}$} &\pzc{$82.837\%$} & \pzc{$5.582\%$}  &  \pzc{$96$}\\ \hline
 \end{tabular}
 \caption{\pzc{Relative offline computational cost of constructing reduced basis and reduced order operators for the lattice problem in Sec. \ref{sec:lattice}. }\label{tab:offline-lattice}}
 \end{table}

\subsubsection{Pin-cell problem \label{sec:pin-cell}}
\begin{figure}[]
  \begin{center} 
\includegraphics[width=0.45\textwidth]{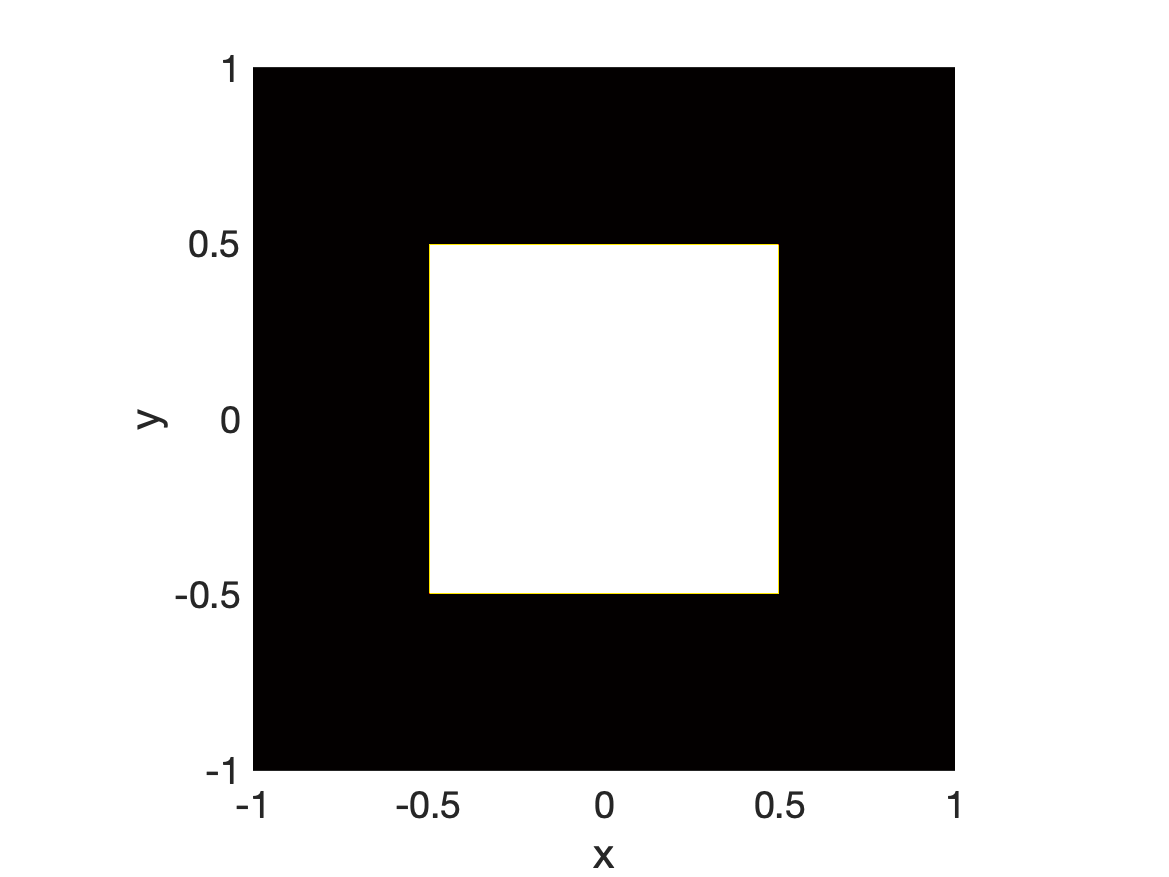}
\includegraphics[width=0.45\textwidth]{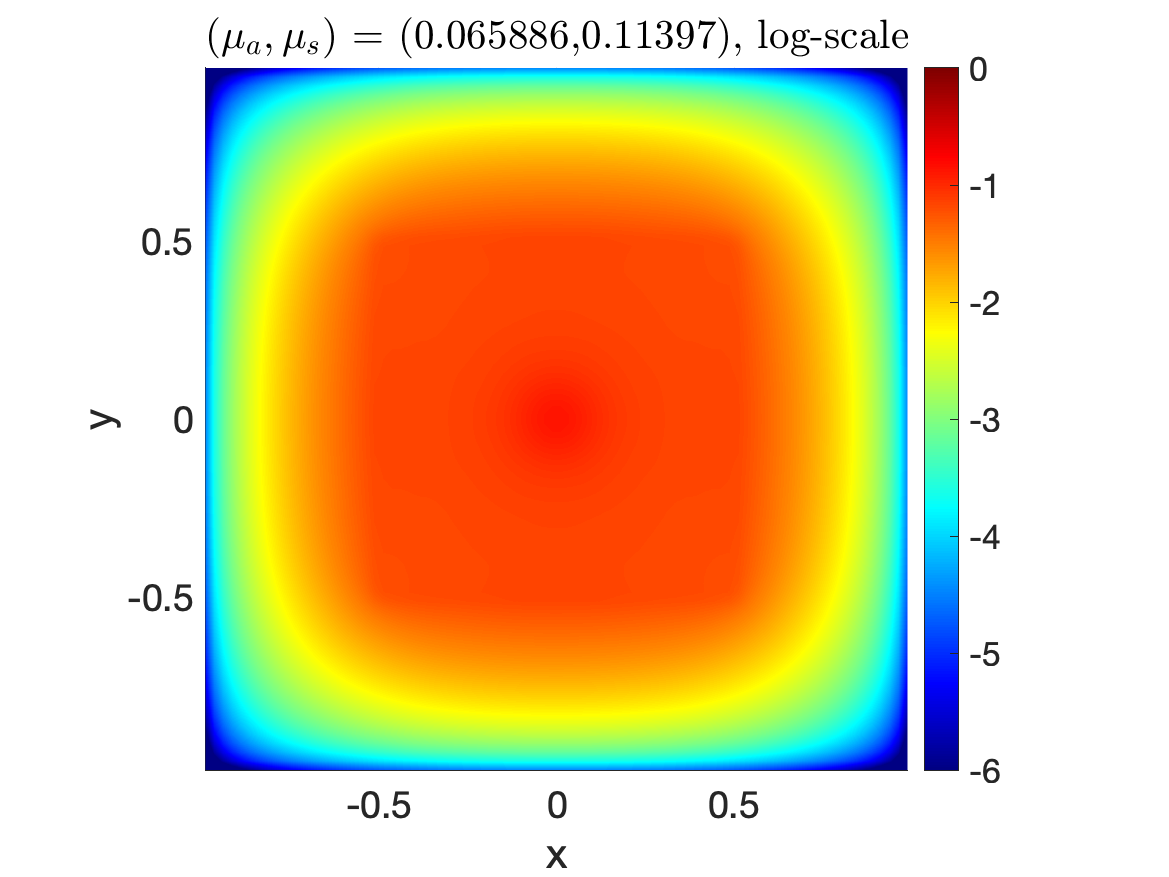}
\includegraphics[width=0.45\textwidth]{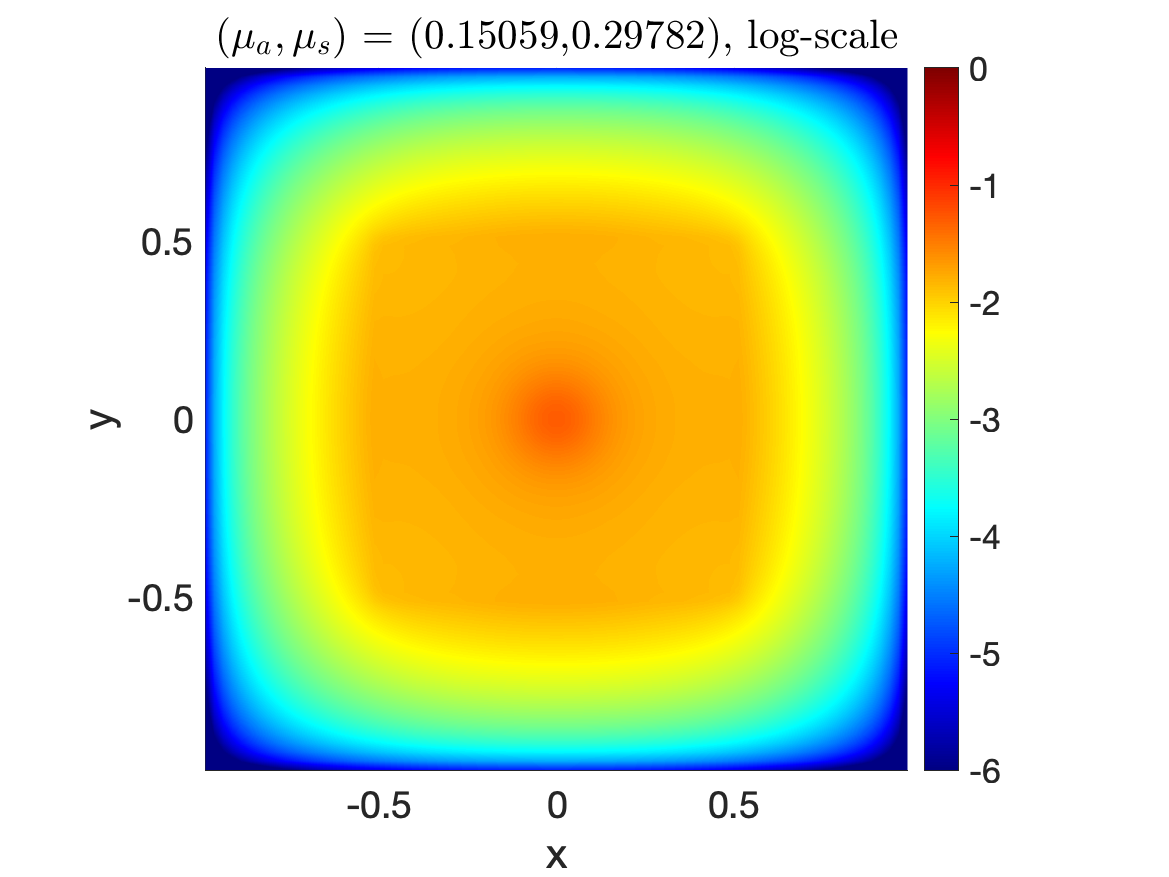}
\includegraphics[width=0.45\textwidth]{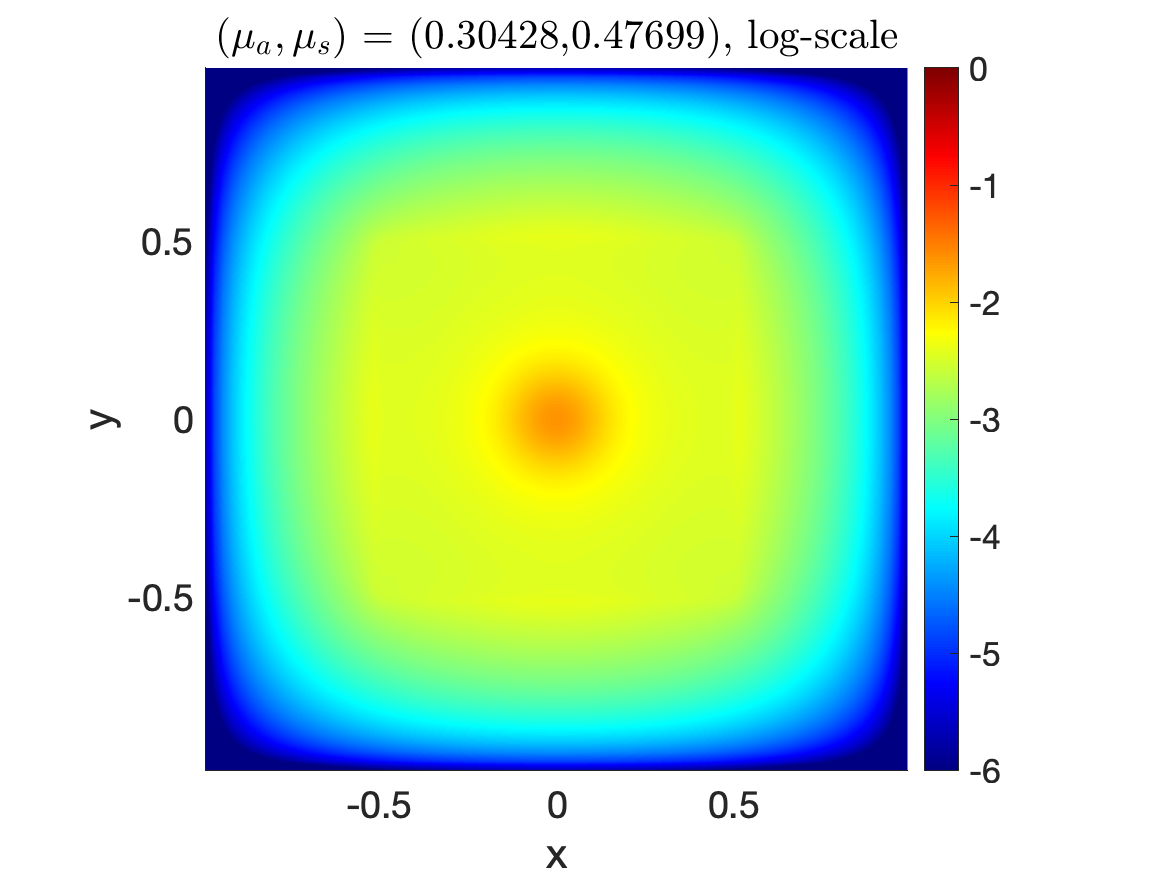}
  \caption{The set-up and reference solutions for the pin-cell problem in Sec. \ref{sec:pin-cell}. Top left: problem set-up. Other pictures: $\rho$ for various pairs of parameters $(\mu_s,\mu_a)$ (under log-scale). \label{fig:pin-cell}}
  \end{center}
\end{figure}

We consider a parametric pin-cell problem with zero inflow boundary conditions and the geometric set-up in the top left picture of Fig. \ref{fig:pin-cell}. The computational domain is $\Omega_{\bx}=[-1,1]^2$ and the source is $G(x,y) = \exp(-100(x^2+y^2))$. The absorption and scattering cross sections are
\begin{equation*}
\sigma_a(x,y) = \begin{cases}
    \mu_a,\;\text{if}\; |x|\leq 0.5\;\text{and}\;|y|\leq0.5,\\
    0,\;\text{otherwise},
\end{cases}\;\text{and}\quad
\sigma_s(x,y) = \begin{cases}
    \mu_s,\;\text{if}\; |x|\leq 0.5\;\text{and}\;|y|\leq0.5,\\
    100,\;\text{otherwise}.
\end{cases}
\end{equation*}
Parameters $\mu_a\in[0.05,0.5]$ and $\mu_s\in[0.05,0.5]$ are the absorption and scattering cross sections for the inner cell, respectively. As we move from the center of the computational domain to the outer part, there is a sharp transition in the strength of the scattering effect from weak to strong.
The computational domain is partitioned with an $80\times80$ uniform mesh. We use the CL($30,6$) quadrature rule for the angular discretization. The tolerance in the stopping criteria for SISA is set to $\epssisa=10^{-11}$. \pzc{The tolerance of relative residual in PGMRES is $2.5\times 10^{-11}$. As shown in Tab. \ref{tab:pin-cell}, converged solutions given by SI and PGMRES result in residuals of comparable magnitude.}

To build ROMs, we use a training set with $25$ uniformly distributed samples: $$\mathcal{P}_{\textrm{train}}=\{(\mu_a,\mu_s)=(i\Delta\mu_a,j\Delta\mu_s),\;i,j=1,\dots 5,\;\Delta\mu_a=\Delta\mu_s=0.05\}.$$ We generate our ROMs with the tolerance  $\epspod=10^{-9}$.
To test the performance of our methods, we randomly sample $10$ pairs of $(\mu_a,\mu_s)$ from $[0.05,0.5]^2$. 
Solutions for $3$ pairs of test samples are presented in Fig. \ref{fig:pin-cell}. When the difference between material properties of the inner and outer parts of the computational domain is larger, the density has sharper features near the material interface. 

\textbf{\pzc{Online efficiency:}} The performance of various methods is demonstrated in Tab. \ref{tab:pin-cell}. 
Due to the significant jump in $\sigma_s$ near the material interface, this problem is challenging for SI-DSA. On average, SI-DSA \pzc{using the zero initial guess} needs $39.3$ iterations to converge. 
\pzc{SI-DSA with ROMIG results in}  approximately $7.6$ times the acceleration compared to DSA and \pzc{nearly $3$ times the acceleration compared to DSA preconditioned GMRES using the zero initial guess.} \pzc{SI-ROMSAD-$3,5$} \pzc{yields} nearly $10$ times the acceleration compared to SI-DSA \pzc{and nearly $4$ times the acceleration compared to DSA preconditioned GMRES using the zero initial guess}.  \pzc{Performance of DSA preconditioned GMRES using ROMIG is comparable with SI-DSA using ROMIG for this example.}

\begin{figure}[]
  \begin{center} 
\includegraphics[width=0.45\textwidth]{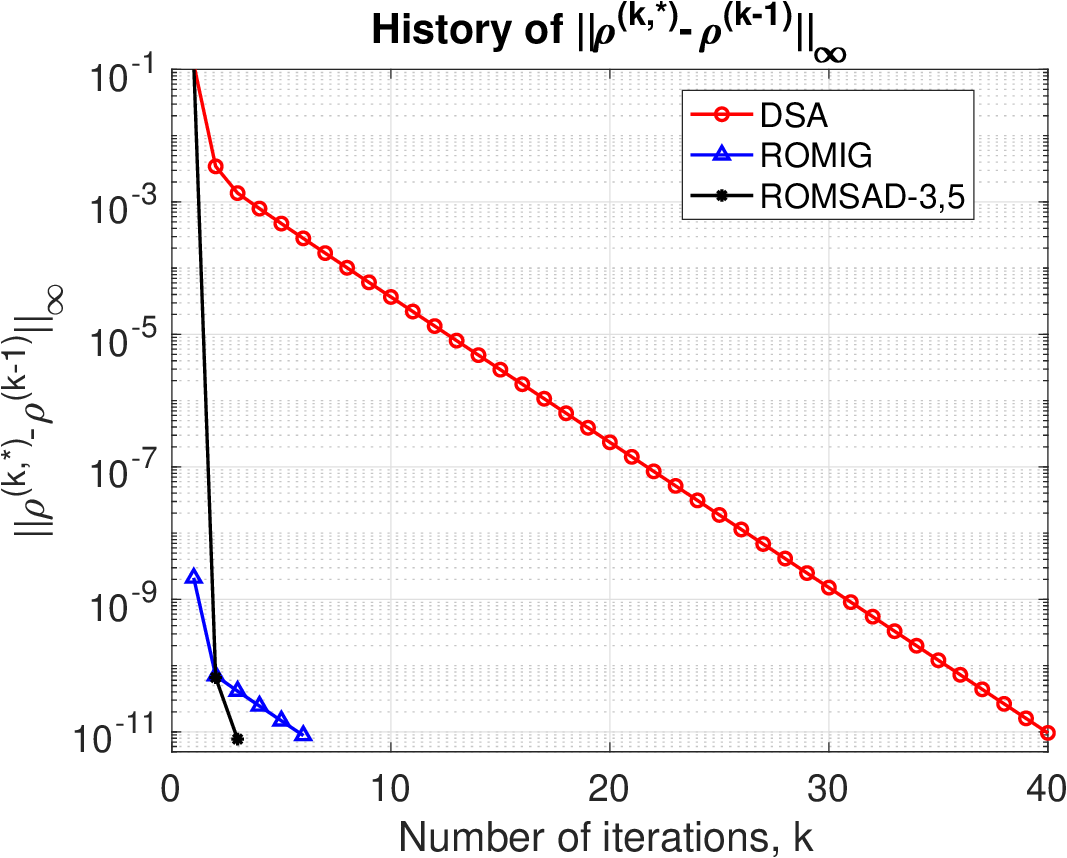}
\includegraphics[width=0.45\textwidth]{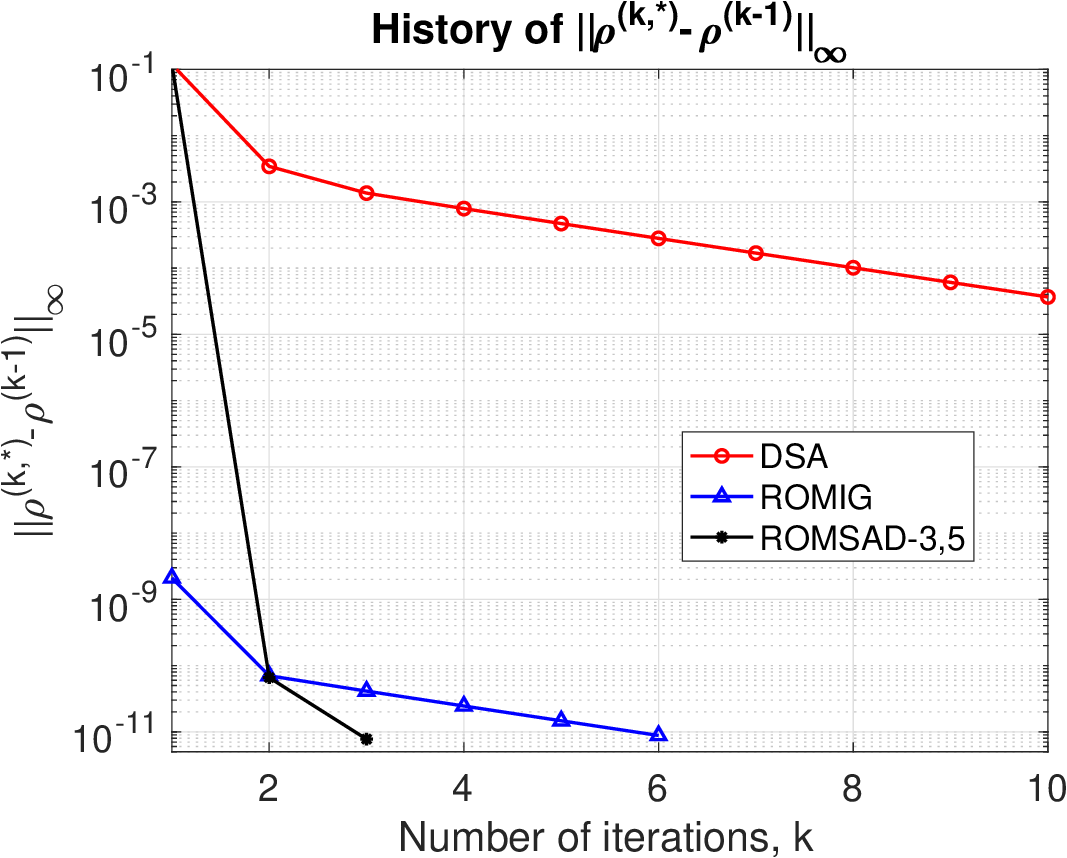}
  \caption{The history of the change of $\brho$ in the $k$-th iteration, i.e. $||\brho^{(k,*)}-\brho^{(k-1)}||_\infty$, for  the pin-cell problem in Sec. \ref{sec:pin-cell} with $(\mu_a,\mu_s)=(0.15068,0.2978)$. Left: entire history. Right: zoomed-in picture. \label{fig:pin_cell_convergence}}
  \end{center}
\end{figure}

In Fig. \ref{fig:pin_cell_convergence}, we show how the change of the density in the $k$-th source iteration, $||\brho^{(k,*)}-\brho^{(k-1)}||_\infty$, evolves for various methods. Due to the use of a ROM-based initial guess, \pzc{ROMIG} starts from a small initial difference $||\brho^{(1,*)}-\brho^{0}||_\infty$, and converges earlier than DSA. The decreasing rates for \pzc{ROMIG} and DSA appear similar to each other. ROMSAD-$3,5$  starts from a zero initial guess, so $||\brho^{(1,*)}-\brho^{(0)}||_\infty$ are the same for ROMSAD-$3,5$ and DSA. Leveraging a reduced order kinetic correction equation, $||\brho^{(k,*)}-\brho^{(k-1)}||_\infty$ for ROMSAD-$3,5$ drops significantly in the second iteration, and its decline rate is faster than DSA in the third iteration. In fact, for the presented test sample, \pzc{with the help of} ROMSAD-$3,5$, \pzc{source iteration} converges with only $3$ iterations.  

\pzc{\textbf{Offline cost:} The main offline costs is solving the RTE for the $25$ training samples.  In Table \ref{tab
}, we present the relative offline computational time for constructing reduced order bases and operators with respect to the average computational time of solving the RTE for one parameter in the test set with SI-DSA. The relative time for constructing reduced basis and operators for the parametric problem itself is only approximately $0.361\%$ and $0.458\%$, respectively. The relative time of constructing reduced basis and operators for the correction equation is only approximately $2.722\%$ and $1.286\%$, respectively.} 

\begin{table}
\centering
 \begin{tabular}{|l|c|c|c|c|c|c|c|c|c|c|c|}
    \hline
    & DSA & \pzc{ROMIG} & \pzc{ROMSAD}-$3,5$ & \pzc{PGMRES} & \pzc{PGMRES-ROMIG}  \\ \hline  
    $\pzc{{n}_{\textrm{sweep}}}$ & 39.3 &  5.3 & 4 & \pzc{16.2} & \pzc{5.6} \\ \hline
    ${T}_{\textrm{rel}}$ & 100\%  & \pzc{13.143\%} & \pzc{9.336\%}   & \pzc{40.070\%} & \pzc{13.781\%} \\ \hline
    $\pzc{\mathcal{R}^\rho_{\infty}}$ & \pzc{2.09e-12} & \pzc{2.23e-12} & \pzc{5.43e-12} & \pzc{8.96e-12} & \pzc{3.99e-12} \\ \hline
 \end{tabular}
 \caption{Results for the pin-cell problem in Sec. \ref{sec:pin-cell}. Dimensions of the reduced order spaces: $r_p=21$ and $r_c=59$. \label{tab:pin-cell}}
 \end{table}
\begin{table}
\centering
 \begin{tabular}{|l|c|c|c|c|c|c|c|c|c|c|c|}
    \hline
    & \pzc{Reduced basis construction} & \pzc{Reduced operator construction}   \\ \hline  
\pzc{Parametric problem} &   \pzc{$0.361\%$} &    \pzc{$0.458\%$}\\ \hline
\pzc{Correction equation} &  \pzc{$2.041\%$} & \pzc{$1.286\%$}  \\ \hline
 \end{tabular}
 \caption{\pzc{Relative offline computational cost of constructing reduced basis and reduced order operators for the pin-cell problem in Sec. \ref{sec:pin-cell}. Dimensions of reduced order spaces: $r_p=21$ and $r_c=59$.}\label{tab:offline-pin-cell}}
 \end{table}

\subsubsection{Variable scattering cross section \label{sec:variable-scattering}}
\begin{figure}[]
  \begin{center} 
\includegraphics[width=0.45\textwidth]{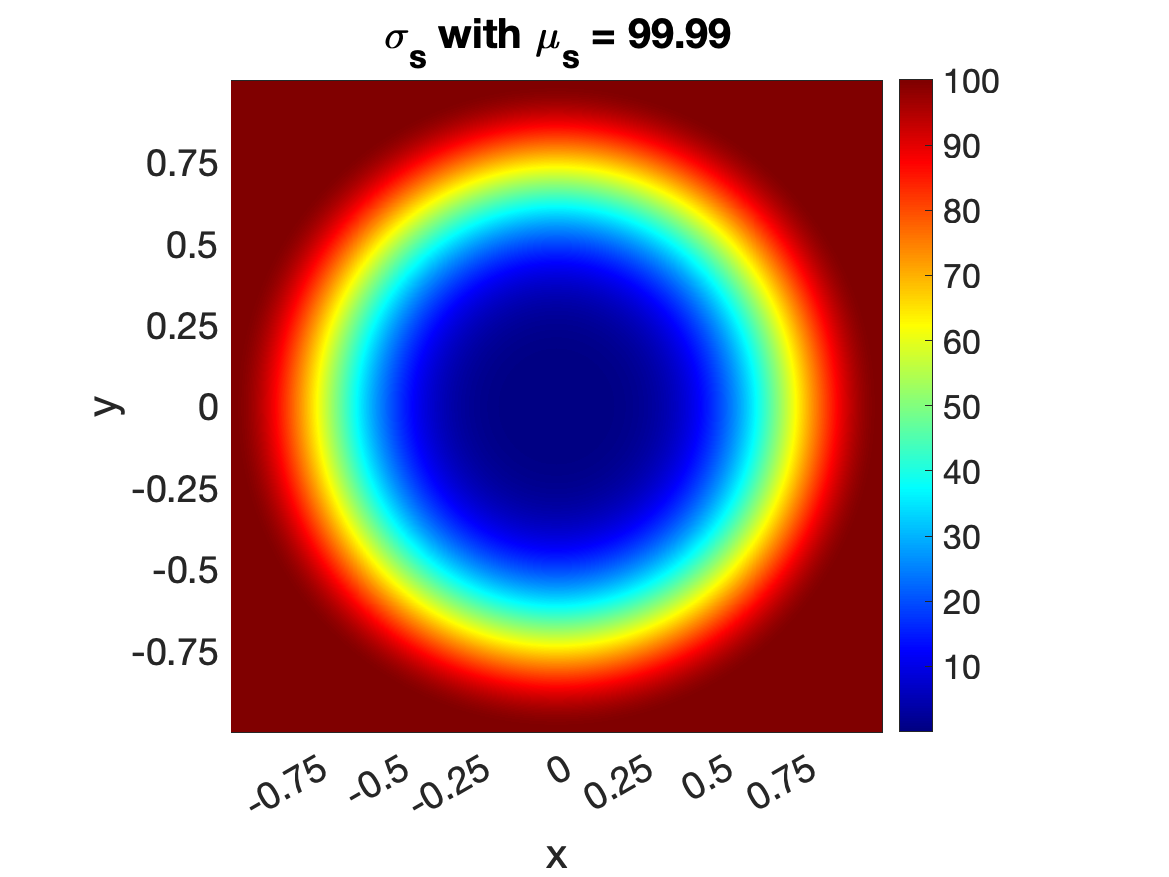}
\includegraphics[width=0.45\textwidth]{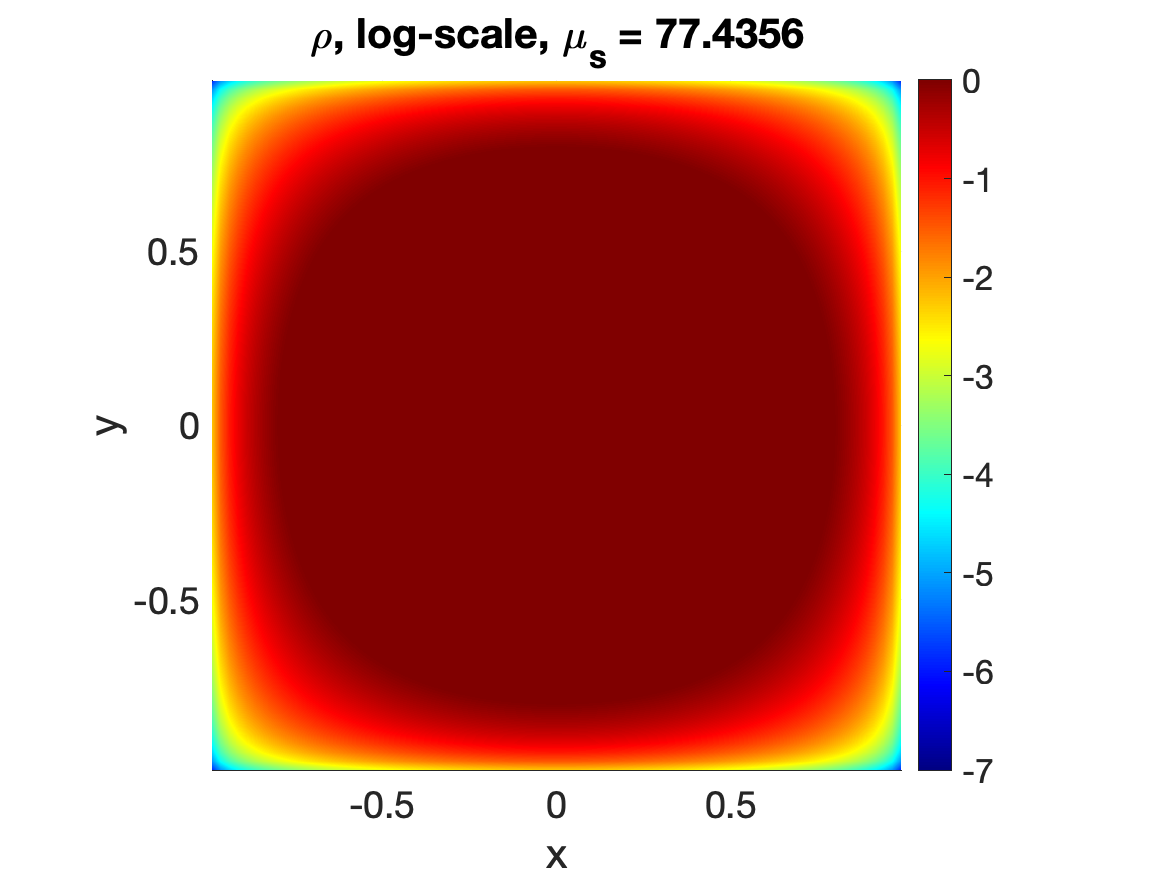}
  \caption{Left: $\sigma_s(\bx)$ with $\mu_s=99.9$ for the variable scattering problem in Sec. \ref{sec:variable-scattering}. Right: $\rho$ for $\mu_s=77.4356$ (under log-scale). \label{fig:variable-scattering}}
  \end{center}
\end{figure}

We consider a parametric problem on the computational domain $\Omega_{\bx}=[-1,1]^2$ with zero inflow boundary conditions and a variable scattering cross section 
\begin{equation}
    \sigma_s(x,y) = \begin{cases}
        \mu_s r^4(2-r^2)^2+0.1,\quad r=\sqrt{x^2+y^2}\leq 1,\\
        \mu_s+0.1,\quad\text{otherwise}, 
    \end{cases}\quad\mu_s\in[49.9,99.9],
\end{equation}
and a zero absorption cross section. The scattering cross section $\sigma_s(\bx)$ with $\mu_s=99.9$ is presented in the left figure of Fig. \ref{fig:variable-scattering}. From the center of the computational domain to the outer part, the scattering cross section smoothly changes from $0.1$ to $\mu_s+0.1$ which is at least $50$ for $\mu_s\in[49.9,99.9]$. In other words, there is a smooth transition from transport dominance to scattering dominance. The parameter $\mu_s$ determines how fast $\sigma_s(\bx)$ changes. The source is $G(x,y)=\frac{10}{\pi}\exp(-100(x^2+y^2)).$ 
We use an $80\times 80$ uniform mesh to partition the computation domain and the CL($30,6$) quadrature rule for the angular discretization. The tolerance in the stopping criteria of SISA is set to $\epssisa=10^{-12}$.
\pzc{The tolerance of relative residual in PGMRES is $2.5\times 10^{-11}$. As shown in Tab. \ref{tab:pin-cell}, converged solutions given by SI and PGMRES result in residuals of comparable magnitude.}

\pzc{\textbf{Online efficiency:}} We use a training set with $50$ uniformly distributed samples and build ROMs with the POD truncation tolerance $\epspod=10^{-11}$. To test the performance of our methods, we randomly sample $10$ values of $\mu_s$ from $[49.9,99.9]$. 
Results for this test are summarized in Tab. \ref{tab:variable-scattering}. Compared with SI-DSA and \pzc{DSA preconditioned GMRES}, all proposed methods result in greater acceleration. \pzc{SI-DSA using ROMIG is the fastest for this example. SI-ROMSAD and DSA preconditioned GMRES using ROMIG have comparable performance. } 

\pzc{\textbf{Offline cost:} The main cost is computing solutions for the $50$ training samples. In Tab. \ref{tab:offline-variable-scattering}, we present the relative offline computational time to construct reduced order basis and operators  with respect to the average computational time of solving RTE for one parameter in the test set with SI-DSA. The relative time for constructing the reduced basis is $3.006\%$ for the parametric problem itself and $8.439\%$ for the correction equation. The relative time for constructing reduced operators is $0.489\%$ for the parametric problem itself and $8.439\%$ for the correction equation.} 
\begin{table}
\centering
 \begin{tabular}{|l|c|c|c|c|c|c|c|c|c|c|c|}
    \hline
    & DSA & \pzc{ROMIG} & \pzc{ROMSAD}-$3,5$ & PGMRES & PGMRES-ROMIG  \\ \hline  
    $\pzc{{n}_{\textrm{sweep}}}$ & 17 &  2.6 & 4.2 & \pzc{12} & \pzc{3.8} \\ \hline
    ${T}_{\textrm{rel}}$ & 100\%   & \pzc{14.357\%} & \pzc{19.879\%} & \pzc{68.760\%} & \pzc{20.627\%}\\ \hline
    $\pzc{\mathcal{R}^\rho_{\infty}}$ & \pzc{2.86e-13} & \pzc{1.14e-13} & \pzc{2.97e-13} & \pzc{7.76e-12} & \pzc{1.12e-12} \\ \hline
 \end{tabular}
 \caption{Results for the variable scattering problem in Sec. \ref{sec:variable-scattering}. Dimensions of the reduced order spaces: $r_p=9$ and $r_c=18$.\label{tab:variable-scattering}}
 \end{table}

\begin{table}
\centering
 \begin{tabular}{|l|c|c|c|c|c|c|c|c|c|c|c|}
    \hline
    & \pzc{Reduced basis construction} & \pzc{Reduced operator construction}   \\ \hline  
\pzc{Parametric problem} &   \pzc{$3.006\%$} &    \pzc{$0.489\%$}\\ \hline
\pzc{Correction equation} &  \pzc{$8.439\%$} & \pzc{$0.703\%$}  \\ \hline
 \end{tabular}
 \caption{\pzc{Relative offline computational cost of constructing reduced basis and reduced order operators for the variable scattering problem in Sec. \ref{sec:variable-scattering}. Dimensions of reduced order spaces: $r_p=9$ and $r_c=18$.}\label{tab:offline-variable-scattering}}
 \end{table}

\section{Conclusion\label{sec:conclusion}}
We propose two strategies to utilize data-driven ROMs to enhance iterative solvers for parametric RTE.
\begin{enumerate}
\item We use the ROM for the parametric problem to provide an improved initial guess for SI-DSA or \pzc{DSA preconditioned Krylov solver}. 
\item We exploit the ROM for the kinetic correction equation to design a new ROM-based synthetic acceleration strategy called ROMSA. We further combine ROMSA with DSA to develop a strategy named ROMSAD. ROMSAD leverages the high efficiency of ROMSA in the early stage of SISA and the robustness of DSA in the later stage. Additionally, we propose an approach to construct the ROM for the  kinetic correction  equation without directly solving it.
\end{enumerate}

In our numerical tests, we \pzc{have the following observations}. 
\begin{enumerate}
    \item \pzc{In the tested examples,}
\pzc{SI-DSA using ROMIG} and \pzc{SI-ROMSAD} \pzc{are more efficient} than SI-DSA or \pzc{DSA preconditioned GMRES  using the zero initial guess.} 
    \item \pzc{In our numerical tests, SI-DSA  with ROMIG is slightly more efficient than or comparable to DSA preconditioned GMRES with ROMIG due to the extra transport sweeps needed in GMRES (see Sec. \ref{sec:details}). However, when the problem  is more challenging or limited data is available, we expect that DSA-preconditioned GMRES using ROMIG will outperform SI-DSA using ROMIG.}
    \item  \pzc{ROMIG} is more robust with respect to different choices of the POD truncation tolerance $\epspod$. ROMSAD has the potential to achieve rapid convergence for some challenging problems.
\end{enumerate} 

Potential future directions are as follows.
\pzc{First, we will investigate how to use ROMs to enhance Krylov solvers, which are preferred for more challenging parametric problem.}
Second, \pzc{the current} switching strategy in \pzc{ROMSAD is ad-hoc to some extent,} \pzc{so we plan to design a more robust strategy} based on \pzc{certified} error estimators. Third, the offline \pzc{efficiency} can be \pzc{improved} by \pzc{leveraging} greedy algorithms \pzc{which adaptively samples training parameters}. Fourth, we plan to extend our method to parametric problems with anisotropic scattering, multiple energy groups and nonlinear terms. Additionally, we aim to integrate our method as a building block for uncertainty quantification, shape optimization, and solvers for inverse problems. 

\pzc{\section*{Acknowledgement}
The author would like to thank Prof. Fengyan Li from the Department of Mathematical Sciences at Rensselaer Polytechnic Institute for her kind suggestions.}

\section*{CRediT authorship contribution statement}

{\bf Zhichao Peng:} Writing – original draft, Writing – review \& editing, Visualization, Validation, Software, Methodology, Data curation, Conceptualization.

\section*{Declaration of generative AI and AI-assisted technologies in the writing process}
During the preparation of this work the author(s) used ChatGPT in order to check grammar errors and improve readability. After using this tool/service, the author(s) reviewed and edited the content as needed and take(s) full responsibility for the content of the publication.

\appendix
\section{\pzc{Algorithm of restarted GMRES}\label{sec:gmres}}

\begin{algorithm}[H]
\caption{\pzc{Outline for preconditioned GMRES with restart every $m$ inner iterations to solve $\BA\bx=\bb$. \cite{barrett1994templates}} \label{alg:gmres} }
\label{alg:source}
\begin{algorithmic}[1]
\STATE{\pzc{Given $\bx^{(0)}$, preconditioner $M^{-1}$, and the maximum number of outer iterations allowed $N_{\textrm{iter}}$.}}
\STATE{\pzc{\textbf{Initialization:} $\br=\BM^{-1}(\bb-\BA\bx)$.}}
\FOR{$l_{\textrm{outer}}=1:N_{\textrm{iter}}$}
    \FOR{$i=1:m$}
        \STATE{\pzc{\textbf{Anorldi process:} 
         (1) Enlarge the Krylov subspace by computing $\bw=\BM^{-1}\BA\bv^{(i)}$. (2) Update the Hessenberg matrix $\BH$ and the orthogonal basis of the Krylov space $\BV$ through modified Gram-Schimidt algorithm.}}
        \STATE{\pzc{Update the residual of the preconditioned system through Givens transformation}. }
        \IF{\pzc{the relative residual is small enough}}
            \STATE{\pzc{\textbf{Update solution through minimization:} $\bx:=\bx^{(0)}+\BV \boldsymbol{y}$, where $\BV \boldsymbol{y}$ is the vector in the Krylov subspace minimizing the residual of the preconditioned system.}}
            \STATE{\pzc{Compute $\br=M^{-1}(\bb-\BA\bx)$.}}
        \ENDIF
    \ENDFOR
    \STATE{\textbf{\pzc{Re-initialization before restart:}} \pzc{update $\bx$ by minimizing the residual and compute $\br=\BM^{-1}(\bb-\BA\bx)$.} }
\ENDFOR
\end{algorithmic}
\end{algorithm}

\section{\pzc{Derivation of consistent DSA}\label{sec:dsa}}
\pzc{Following \cite{adams2001discontinuous}, we briefly outline how the consistent discretization of DSA is derived. For simplicity, we focus on the 1D slab geometry. We also assume that if $v_j$ is a quadrature point, then $-v_j$ is also a quadrature point with the same weight.}

\pzc{Let $\{T_i=[x_{i-\half},x_{i+\half}], \; i=1,\dots, N_x\}$ be a partition of the computational domain. Let $\{\phi(\bx)\}_{i=1}^{N_{\textrm{DOF}}}$ be an orthonormal basis for the polynomial space used by the upwind DG spatial discretization. Define discrete operators $\BD^\pm \in\mathbb{R}^{N_{\textrm{DOF}}\times N_{\textrm{DOF}}}$ as follows:
\begin{subequations}
\begin{align}
(\BD^+)_{kl}&=-\sum_{i=1}^{N_x}\int_{T_i} \partial_x\phi_k(x) \phi_l(x) dx +\sum_{i=1}^{N_x-1}\phi_l(x^+_{i+\half}) \phi_k(x^-_{i+\half})-\sum_{i=1}^{N_x}\phi_l(x^+_{i-\half}) \phi_k(x^+_{i-\half}), \\
(\BD^-)_{kl}&=-\sum_{i=1}^{N_x}\int_{T_i} \partial_x\phi_k(x) \phi_l(x) dx +\sum_{i=1}^{N_x}\phi_l(x^-_{i+\half}) \phi_k(x^-_{i+\half})-\sum_{i=2}^{N_x}\phi_l(x^-_{i-\half}) \phi_k(x^+_{i-\half}).
\end{align}
\end{subequations}
Let $\BD_C=\frac{1}{2}(\BD^++\BD^-)$ and $\BD_J=\BD^+-\BD^-$. Then, the discrete operator for the upwind DG discretization can be represented as
$$ \BD_j = v_j\left(\BD_C-\frac{\textrm{sign}(v_j)}{2}D_J\right).$$
}

\pzc{
Assume $\delta f^{(l)}(x,v_j) = \delta\rho^{(l)}(x)+3v_j \delta J^{(l)}(x)$. The matrix-vector form of the upwind DG discretization of the kinetic correction equation becomes
\begin{subequations}
\label{eq:upwind_dg_dsa_correction}
\begin{align}
&\left(v_j (\BD_C- \frac{1}{2}\BD_J)+\BSigma_t\right)\left(\delta\brho^{(l)}+3v_j\delta\BJ^{(l)}\right) = \BSigma_s \delta\brho + \BSigma_s (\brho^{(l,*)}-\brho^{(l-1)}),\quad \text{if}\; v_j\geq 0,\\
&\left(v_j (\BD_C+\frac{1}{2}\BD_J)+\BSigma_t\right)\left(\delta\brho^{(l)}+3v_j\delta\BJ^{(l)}\right) = \BSigma_s \delta\brho + \BSigma_s (\brho^{(l,*)}-\brho^{(l-1)}),\quad \text{if}\; v_j<0.
\end{align}
\end{subequations}
We then numerically compute the zero-th and first order moment of \eqref{eq:upwind_dg_dsa_correction} in the angular space: 
\begin{subequations}
\label{eq:dsa_first_moment}
\begin{align}
    &\BD_C\delta\BJ^{(l)} - (\sum_{v_j>0} \omega_jv_j) \BD_J\delta\brho^{(l)}+\BSigma_a \delta\brho^{(l)} = \BSigma_s (\brho^{(l,*)}-\brho^{(l-1)}),\label{eq:dsa_rho}\\
    (\sum_{j=1}^{N_v}\omega_j v_j^2)&\BD_C\delta\brho^{(l)} +\left(\BSigma_t- 3(\sum_{v_j>0}\omega_j v_j^3)\BD_J\right)\delta\BJ^{(l)} = 0.\label{eq:dsa_J}
\end{align}
\end{subequations}
Equation \eqref{eq:dsa_J} yields
\begin{equation}
\delta\BJ^{(l)}=-(\sum_{j=1}^{N_v}\omega_j v_j^2)\left(\BSigma_t- 3(\sum_{v_j>0}\omega_j v_j^3)\BD_J\right)^{-1}\BD_C\delta\brho^{(l)}.
\label{eq:dsa_J_in_rho}
\end{equation} 
The consistent discretization for the diffusion correction equation \eqref{eq:diffusion_correction} equation can be obtained by substituting \eqref{eq:dsa_J_in_rho} into \eqref{eq:dsa_rho}.
}

\bibliographystyle{elsarticle-num} 
\bibliography{ref}
\end{document}